\documentclass[11pt]{article}
\usepackage{amsfonts}
\usepackage{epsfig}
\usepackage{color}
\def\be{\begin{equation}}
\def\ee{\end{equation}}
\def\bea{\begin{eqnarray}}
\def\eea{\end{eqnarray}}
\def\bes{\begin{eqnarray*}}
\def\ees{\end{eqnarray*}}
\def\pmatrix{\left(\begin{array}{cc}}
\def\endpmatrix{\end{array}\right)}

\def\sgn{{\rm sgn}}

\def\Gr{{\rm Gr}}
\def\Lag {{\rm Lag}}
\def\Gr{{\rm Gr}}
\def\sign{{\rm sign}}
\def\det{{\rm det}}

\def\exp {{\rm exp}}

\def\Sp{{\rm Sp}}

\def\rank{{\rm rank}}

\def\nn{\nonumber}
\def\<{\langle}
\def\>{\rangle}
\def\lb{\label}
\def\bs{\setminus}

\def\R{{\bf R}}
\def\C{{\bf C}}
\def\Z{{\bf Z}}
\def\N{{\bf N}}
\def\U{{\bf U}}
\def\Q{{\bf Q}}

\def\ga{{\gamma}}

\def\om{{\omega}}
\def\Om{{\Omega}}

\def\var{{\varepsilon}}
\def\lm{{\lambda}}
\def\Lm{{\Lambda}}

\def\sg{{\sigma}}
\def\Sg{{\Sigma}}

\def\G{{\cal G}}

\def\P{{\cal P}}
\def\J{{\cal J}}

\def\I{{\cal I}}

\def\r{\color{red}}
\def\td#1{\tilde{#1}}
\def\hb{\vrule height0.18cm width0.14cm $\,$}

\def\td#1{\tilde{#1}}
\def\diag{{\rm diag}}

\title{Multiple brake orbits on compact convex
symmetric reversible hypersurfaces in $\R^{2n}$ }

\author{Duanzhi Zhang
\thanks{Partially supported by the NSF of China
(10801078, 11171314) and Nankai University. E-mail:
zhangdz@nankai.edu.cn}\qquad and \qquad Chungen Liu
\thanks{Corresponding author. Partially supported by the NSF of China (11071127, 10621101), 973 Program of
MOST (2011CB808002). E-mail:
liucg@nankai.edu.cn}\\ \\
School of Mathematics and LPMC, Nankai University\\
Tianjin 300071, People's Republic of China}

\date{}
\begin{document}
\maketitle  \begin{abstract} In this paper, we prove that there
exist at least $\left[\frac{n+1}{2}\right]+1$ geometrically distinct
brake orbits on every $C^2$ compact convex symmetric hypersurface
$\Sg$ in $\R^{2n}$  for $n\ge 2$ satisfying the reversible condition
$N\Sg=\Sg$ with $N=\diag (-I_n,I_n)$. As a consequence, we show that
there exist at least $\left[\frac{n+1}{2}\right]+1$ geometrically
distinct brake orbits in every bounded convex symmetric domain in
$\R^{n}$ with $n\ge 2$ which gives a positive answer to the Seifert
conjecture of 1948 in the symmetric case for $n=3$. As an
application, for $n=4$ and $5$, we prove that if there are exactly
$n$ geometrically distinct closed characteristics on $\Sg$, then all
of them are symmetric brake orbits after suitable time translation.
\end {abstract}
\noindent {\bf Key
words:} { Brake orbit, Maslov-type index, H${\rm \ddot{o}}$rmander
index, Convex symmetric}

\renewcommand{\theequation}{\thesection.\arabic{equation}}

\setcounter{equation}{0}
\section{Introduction }

 Let $V\in C^2(\R^n, \R)$ and $h>0$ such that $\Om\equiv \{q\in
\R^n|V(q)<h\}$ is nonempty, bounded, open and connected. Consider
the following fixed energy problem of the second order autonomous
Hamiltonian system
 \bea && \ddot{q}(t)+V'(q(t))=0, \quad {\rm for}\;q(t)\in\Om, \lb{1.1}\\
&& \frac{1}{2}|\dot{q}(t)|^2+V(q(t))= h, \qquad\forall t\in\R, \lb{1.2}\\
&& \dot{q}(0)=\dot{q}(\frac{\tau}{2})=0, \lb{1.3}\\
&& q(\frac{\tau}{2}+t)=q(\frac{\tau}{2}-t),\qquad q(t+\tau)=q(t),
    \quad \forall t\in\R.  \lb{1.4}
\eea

 A solution $(\tau,q)$ of (\ref{1.1})-(\ref{1.4}) is called a {\it
brake orbit$\;$} in $\Om$. We call two brake orbits $q_1$ and
$q_2:\R\to\R^n$ {\it geometrically distinct}  if $q_1(\R)\neq
q_2(\R)$.

We denote by $\mathcal{O}(\Om)$ and $\td{\mathcal{O}}(\Om)$ the sets
of all brake orbits and geometrically distinct brake orbits in $\Om$
respectively.

Let $J_k=\left(\begin{array}{cc}0&-I_k\\I_k&0\end{array}\right)$ and
$N_k=\left(\begin{array}{cc}-I_k&0\\0&I_k\end{array}\right)$ with
$I_k$ being the identity in $\R^k$. If $k=n$ we will omit the
subscript $k$ for convenience, i.e., $J_n=J$ and $N_n=N$.

 The symplectic group $\Sp(2k)$ for any $k\in\N$ is defined by
 $$\Sp(2n)=\{M\in \mathcal{L}(\R^{2k})| M^TJ_kM=J_k\},$$
where $M^T$ is the transpose of matrix $M$.

For any $\tau>0$, the symplectic path in $\Sp(2k)$ starting from the
identity $I_{2k}$ is defined by
  $$\mathcal{P}_{\tau}(2k)=\{\gamma\in C([0,\tau],\Sp(2k))|
  \gamma(0)=I_{2k}\}.$$

  Suppose that $H\in
C^2(\R^{2n}\bs\{0\},\R)\cap C^1(\R^{2n},\R)$ satisfying

\be H(Nx)=H(x),\qquad \forall\, x\in\R^{2n}.\lb{1.5}\ee

We consider the following fixed energy problem
 \bea
\dot{x}(t) &=& JH'(x(t)), \lb{1.6}\\
H(x(t)) &=& h,   \lb{1.7}\\
x(-t) &=& Nx(t),  \lb{1.8}\\
 x(\tau+t) &=& x(t),\; \forall\,t\in\R. \lb{1.9} \eea

A solution $(\tau,x)$ of (\ref{1.6})-(\ref{1.9}) is also called a
{\it brake orbit} on  $\Sg:=\{y\in\R^{2n}\,|\, H(y)=h\}$.

\noindent{\bf Remark 1.1.} It is well known that via \be H(p,q) =
{1\over 2}|p|^2 + V(q), \lb{1.10}\ee  $x=(p,q)$ and $p=\dot q$, the
elements in $\mathcal{O}(\{V<h\})$ and the solutions of
(\ref{1.6})-(\ref{1.9}) are one to one correspondent.

In more general setting, let $\Sg$ be a $C^2$ compact hypersurface
in $\R^{2n}$ bounding a compact set $C$ with nonempty interior.
Suppose $\Sg$ has non-vanishing Guassian curvature  and satisfies
the reversible condition $N(\Sg-x_0)=\Sg-x_0:=\{x-x_0|x\in \Sg\}$
for some $x_0\in C$. Without loss of generality, we may assume
$x_0=0$. We denote the set of all such hypersurface in $\R^{2n}$ by
$\mathcal{H}_b(2n)$. For $x\in \Sg$, let $N_\Sg(x)$ be the unit
outward
 normal  vector at $x\in \Sg$. Note that here by the reversible
condition there holds $N_\Sg(Nx)=NN_\Sg(x)$. We consider the
dynamics problem of finding $\tau>0$ and an absolutely continuous
curve $x:[0,\tau]\to \R^{2n}$ such that
  \bea  \dot{x}(t)&=&JN_\Sg(x(t)), \qquad x(t)\in \Sg,\lb{1.11}\\
      x(-t) &=& Nx(t), \qquad x(\tau+t) = x(t),\qquad {\rm
  for\;\; all}\;\; t\in \R.\lb{1.12}\eea

A solution $(\tau,x)$ of the problem (\ref{1.11})-(\ref{1.12}) is a
special closed characteristic on $\Sg$, here we still call it a
brake orbit on $\Sg$.

We also call two brake orbits $(\tau_1, x_1)$ and $(\tau_2,x_2)$
{\it geometrically distinct} if $x_1(\R)\ne x_2(\R)$, otherwise we
say they are equivalent. Any two equivalent brake orbits are
geometrically the same. We denote by ${\mathcal{J}}_b(\Sg)$ the set
of all brake orbits on $\Sg$, by $[(\tau,x)]$ the equivalent class
of $(\tau,x)\in {\mathcal{J}}_b(\Sg)$ in this equivalent relation
and by $\td{\mathcal{J}}_b(\Sg)$ the set of $[(\tau,x)]$ for all
$(\tau,x)\in {\mathcal{J}}_b(\Sg)$. From now on, in the notation
$[(\tau,x)]$ we always assume $x$ has minimal period $\tau$. We also
denote by $\tilde {\mathcal {J}}(\Sg)$ the set of all geometrically
distinct closed characteristics on $\Sg$.

Let $(\tau,x)$ be a solution of (\ref{1.6})-(\ref{1.9}). We consider
the boundary value problem of the linearized Hamiltonian system \bea
&&\dot{y}(t) = JH''(x(t))y(t),  \lb{1.13}\\
&&y(t+\tau)=y(t), \quad y(-t)=Ny(t), \qquad \forall t\in\R.
\lb{1.14} \eea

Denote by $\ga_x(t)$ the fundamental solution of the system
(\ref{1.13}), i.e., $\ga_x(t)$ is the solution of the following
problem
 \bea
\dot{\ga_x}(t) &=& JH''(x(t))\ga_x(t), \lb{1.15}\\
\ga_x(0) &=& I_{2n}.  \lb{1.16} \eea
 We call $\ga_x\in C([0,\tau/2],\Sp(2n))$ the {\it
associated symplectic path} of $(\tau, x)$.

Let $B^n_1(0)$ denote the open unit ball $\R^n$ centered at the
origin $0$. In \cite{Se1} of 1948, H. Seifert proved
$\td{\mathcal{O}}(\Om)\neq \emptyset$ provided $V'\neq 0$ on
$\partial \Om$, $V$ is analytic and $\Om$ is homeomorphic to
$B^n_1(0)$. Then he proposed his famous conjecture: {\it
$^{\#}\tilde{\mathcal{O}}(\Om)\geq n$  under the same conditions}.

After 1948, many studies have been carried out for the brake orbit
problem. S. Bolotin proved first in \cite{Bol}(also see \cite{BolZ})
of 1978 the existence of brake orbits in general setting.
 K. Hayashi in \cite {Ha1}, H. Gluck and W. Ziller
 in \cite{GZ1}, and V. Benci in \cite {Be1} in 1983-1984
  proved $^{\#}\td{\mathcal{O}}(\Om)\geq 1$ if $V$ is
$C^1$, $\bar{\Om}=\{V\leq h\}$ is compact, and $V'(q)\neq 0$ for all
$q\in \partial{\Om}$. In 1987, P. Rabinowitz in \cite{Ra1} proved
that if $H$ satisfies (\ref{1.5}), $\Sg\equiv H^{-1}(h)$ is
star-shaped, and $x\cdot H'(x)\neq 0$ for all $x\in \Sg$, then
$^{\#}\td{\mathcal{J}}_b(\Sg)\geq 1$. In 1987, V. Benci and F.
Giannoni gave a different proof of the existence of one brake orbit
in \cite{BG}.

In 1989, A. Szulkin in \cite{Sz} proved that
$^{\#}\td{\J_b}(H^{-1}(h))\geq n$, if $H$ satisfies conditions in
\cite{Ra1} of Rabinowitz and the energy hypersurface $H^{-1}(h)$ is
$\sqrt{2}$-pinched. E. van Groesen in \cite{Gro} of 1985 and A.
Ambrosetti, V. Benci, Y. Long in \cite{ABL1} of 1993 also proved
$^{\#}\td{\mathcal{O}}(\Om)\geq n$ under different pinching
conditions.

Without pinching condition, in \cite{LZZ} Y. Long, C. Zhu and the
second author of this paper proved the following result:
 {\it For $n\ge 2$, suppose $H$ satisfies

(H1) (smoothness) $H\in C^2(\R^{2n}\bs\{0\},\R)\cap
C^1(\R^{2n},\R)$,

(H2) (reversibility) $H(Ny)=H(y)$ for all $y\in\R^{2n}$.

(H3) (convexity) $H''(y)$ is positive definite for all
$y\in\R^{2n}\bs\{0\}$,

(H4) (symmetry) $H(-y)=H(y)$ for all $y\in\R^{2n}$.

\noindent Then for any given $h>\min \{ H(y)|\; y\in \R^{2n}\}$ and
$\Sg=H^{-1}(h)$, there holds } $$^{\#}\td{\J}_b(\Sg)\ge 2.$$

As a consequence they also proved that: {\it For $n\geq 2$, suppose
$V(0)=0$, $V(q)\geq 0$, $V(-q)=V(q)$ and $V''(q)$ is positive
definite for all $q\in \R^n\bs\{0\}$. Then for $\Om\equiv
\{q\in\R^n|V(q)<h\}$ with $h>0$, there holds}
$$^{\#}\td{\mathcal{O}}(\Om)\ge 2.$$

Under the same condition of \cite{LZZ}, in 2009 Liu and Zhang in
\cite{LiuZhang} proved that $^{\#}\td{\J}_b(\Sg)\ge
\left[\frac{n}{2}\right]+1$, also they proved
$^{\#}\td{\mathcal{O}}(\Om)\ge \left[\frac{n}{2}\right]+1$ under the
same condition of \cite{LZZ}. Moreover if all brake orbits on $\Sg$
are nondegenerate, Liu and Zhang in \cite{LiuZhang} proved that
         $^{\#}\td{\J}_b(\Sg)\ge n+\mathfrak{A}({\Sg}),$
where $2\mathfrak{A}(\Sigma)$ is the number of  geometrically
distinct asymmetric brake orbits on $\Sg$.

\noindent{\bf Definition 1.1.} {\it We denote
$$\begin{array}{ll}\mathcal{H}_b^{c}(2n)=\{\Sg\in \mathcal{H}_{b}(2n)|\;\Sg\; { is\;
strictly\; convex\;} \},\\\mathcal{H}_b^{s,c}(2n)=\{\Sg\in
\mathcal{H}_{b}^c(2n)| \;-\Sg=\Sg\}.\end{array}$$}

\noindent{\bf Definition 1.2.} {\it For
$\Sg\in\mathcal{H}_b^{s,c}(2n)$, a brake orbit $(\tau,x)$ on $\Sg$
 is called symmetric if
$x(\R)=-x(\R)$. Similarly, for a $C^2$ convex symmetric bounded
domain $\Omega\subset \R^n$, a brake orbit $(\tau,q)\in \mathcal
{O}(\Omega)$ is called symmetric if $q(\R)=-q(\R)$.}

Note that a brake orbit $(\tau,x)\in \mathcal {J}_b(\Sg)$ with
minimal period $\tau$
 is symmetric if
$x(t+\tau/2)=-x(t)$ for $t\in \R$, a brake orbit $(\tau,q)\in
\mathcal {O}(\Omega)$ with minimal period $\tau$ is  symmetric if
$q(t+\tau/2)=-q(t)$ for $t\in\R$.

In this paper, we denote by $\N$, $\Z$, $\Q$ and $\R$ the sets of
positive integers, integers, rational numbers and real numbers
respectively. We denote by $\langle \cdot,\cdot\rangle$ the standard
inner product in $\R^n$ or $\R^{2n}$, by $(\cdot,\cdot)$ the inner
product of corresponding Hilbert space. For any $a\in \R$, we denote
$E(a)=\inf\{k\in \Z|k\ge a\}$ and $[a]=\sup\{k\in \Z|k\le a\}$.

The following are the main results of this paper.

\noindent{\bf Theorem 1.1.} {\it For any
$\Sg\in\mathcal{H}_b^{s,c}(2n)$ with $n\ge 2$, we have }
           $$^{\#}\td{\J}_b(\Sg)\ge \left[\frac{n+1}{2}\right]+1 .$$

\noindent{\bf Corollary 1.1.} {\it Suppose $V(0)=0$, $V(q)\geq 0$,
$V(-q)=V(q)$ and $V''(q)$ is positive definite for all $q\in
\R^n\bs\{0\}$ with $n\ge 3$. Then for any given $h>0$ and $\Om\equiv
\{q\in\R^n|V(q)<h\}$, we
  have}
$$^{\#}\td{\mathcal{O}}(\Om)\ge \left[\frac{n+1}{2}\right]+1.$$

\noindent{\bf Remark 1.2.} Note that for $n=3$, Corollary 1.1 yields
$^{\#}\td{\mathcal{O}}(\Om)\ge 3$, which gives a positive answer to
Seifert's conjecture in the convex symmetric case.

As a consequence of Theorem 1.1, we can prove

\noindent{\bf Theorem 1.2.} {\it For $n=4,5$ and any
$\Sg\in\mathcal{H}_b^{s,c}(2n)$, suppose
         $$^{\#}\td{\J}(\Sg)= n.$$
Then all of them are symmetric brake orbits after suitable
translation.}

\noindent{\bf Example 1.1.} A typical example of $\Sg\in
\mathcal{H}_b^{s,c}(2n)$ is the ellipsoid $\mathcal {E}_n(r)$
defined as follows. Let $r=(r_1,\cdots,r_n)$ with $r_j>0$ for $1\le
j\le n$. Define
$$\mathcal {E}_n(r)=\left\{x=(x_1,\cdots,x_n, y_1,\cdots,y_n)\in\R^{2n}\;
\left|\;\sum_{k=1}^n\frac{x_k^2+y_k^2}{r_k^2}=1\right.\right\}.$$
 If $r_j/r_k\notin \Q$ whenever $j\ne k$, from \cite{Ek} one can see
 that there are precisely $n$ geometrically distinct symmetric brake orbits on $\mathcal
 {E}_n(r)$ and all of them are nondegenerate.

\setcounter{equation}{0}
\section {Index theories of $(i_{L_j},\nu_{L_j})$ and $(i_\om,\nu_\om)$}  

Let $ \mathcal {L}(\R^{2n})$ denotes the set of $2n\times 2n$  real
matrices  and $\mathcal{L}_s(\R^{2n})$ denotes its subset of
symmetric ones. For any $F\in \mathcal{L}_s(\R^{2n})$, we denote by
$m^*(F)$ the dimension of maximal positive
     definite subspace, negative definite subspace, and kernel of
     any $F$ for $*=+,-,0$ respectively.

In this section, we make some preparation for the proof of Theorem
3.1 below. We first briefly review the index function
 $(i_\om,\nu_\om)$ and $(i_{L_j},\nu_{L_j})$ for $j=0,1$, more details can be found in \cite{LiuZhang}
and \cite{Long1}. Following Theorem 2.3 of \cite{Zhang1} we study
the differences $i_{L_0}(\ga)-i_{L_1}(\ga)$ and
$i_{L_0}(\ga)+\nu_{L_0}(\ga)-i_{L_1}(\ga)-\nu_{L_1}(\ga)$ for
$\ga\in\mathcal{P}_\tau(2n)$ by compute $\sgn M_\var(\ga(\tau))$. We
obtain some basic lemmas which will be used frequently in the proof
of the main theorem of this paper.

 For any
$\om\in \U$, the following codimension 1 hypersuface in $\Sp(2n)$ is
defined by:
    $$ \Sp(2n)_\om^0=\{M\in \Sp(2n)|\det(M-\om I_{2n})=0\}.$$
For any two continuous path $\xi$ and $\eta$: $[0,\tau]\to\Sp(2n)$
with $\xi(\tau)=\eta(0)$, their joint path is defined by
    \bea  \eta * \xi(t)=\left\{\begin{array}{lr} \xi(2t)\qquad &{\rm
     if}\,
     0\le t\le \frac{\tau}{2},\\ \eta(2t-\tau)\quad &{\rm if}\,
     \frac{\tau}{2}\le t\le \tau .\end{array}\right.\nn\eea
Given any two $(2m_k\times 2m_k)$- matrices of square block form
$M_k=\left(\begin{array}{cc}A_k&B_k\\C_k&D_k\end{array}\right)$ for
$k=1,2$, as in \cite{Long1}, the $\diamond$-product of $M_1$ and
$M_2$ is defined by the following $(2(m_1+m_2)\times
2(m_1+m_2))$-matrix $M_1\diamond M_2$:
     $$M_1\diamond M_2=\left(\begin{array}{cccc}A_1&0&B_1&0\\0&A_2&0&B_2\\
        C_1&0&D_1&0\\0&C_2&0&D_2\end{array}\right).$$
A special path $\xi_n$ is defined by
        $$\xi_n(t)=\left(\begin{array}{cc}2-\frac{t}{\tau}&0\\0&(2-\frac{t}{\tau})^{-1}\end{array}\right)^{\diamond n}, \qquad \forall t\in[0,\tau].$$
{\bf Definition 2.1.} For any $\om\in\U$ and $M\in\Sp(2n)$, define
     \bea \nu_\om(M)=\dim_\C\ker(M-\om I_{2n}).\nn\eea
For any $\ga\in \mathcal{P}_\tau(2n)$, define
  \bea \nu_\om(\ga)=\nu_\om(\ga(\tau)).\nn\eea
If $\ga(\tau)\notin\Sp(2n)_\om^0$, we define
     \be i_\om(\ga)=[\Sp(2n)_\om^0\,:\,\ga*\xi_n],\lb{m1}\ee
where the right-hand side of (\ref{m1}) is the usual homotopy
intersection number and the orientation of $\ga*\xi_n$ is its
positive time direction under homotopy with fixed endpoints. If
$\ga(\tau)\in\Sp(2n)_\om^0$, we let $\mathcal{F}(\ga)$ be the set of
all open neighborhoods of $\ga$ in $\mathcal{P}_\tau(2n)$, and
define
 \bea i_\om(\ga)=\sup_{U\in\mathcal{F}(\ga)}\inf\{i_\om(\beta)|\,\beta(\tau)\in U \,{\rm and}\,
 \beta(\tau)\notin\Sp(2n)_\om^0\}.\nn\eea
Then $(i_\om(\ga),\nu_\om(\ga))\in \Z\times \{0,1,...,2n\}$, is
called the index function of $\ga$ at $\om$.

For any $M\in\Sp(2n)$ we define
 \bea \Om(M)=\{P\in \Sp(2n)&|&\sg(P)\cap\U=\sg(M)\cap \U\nn\\
                            && {\rm and}\,
                            \nu_\lm(P)=\nu_\lm(M),\;\; \forall
                            \lm\in\sg(M)\cap\U\},\nn\eea
where we denote by $\sg(P)$ the spectrum of $P$.

 We denote by $\Om^0(M)$ the path connected component of
$\Om(M)$ containing $M$, and call it the {\it homotopy component} of
$M$ in $\Sp(2n)$.

\noindent{\bf Definition 2.2.} For any $M_1$,$M_2\in \Sp(2n)$, we
call $M_1\approx M_2$ if $M_1\in \Om^0(M_2)$.

\noindent{\bf Remark 2.1.} It is easy to check that $\approx$ is an
equivalent relation. If $M_1\approx M_2$, we have $M_1^k\approx
M_2^k$ for any $k\in \N$ and $M_1\diamond M_3 \approx M_2\diamond
M_4$ for $M_3\approx M_4$. Also we have $PMP^{-1}\approx M$ for any
$P,M\in \Sp(2n)$.

The following symplectic matrices were introduced as {\it basic
normal forms} in \cite{Long1}:
      \bea D(\lm)=
      \left(\begin{array}{cc}\lm&0\\0&\lm^{-1}\end{array}\right),\qquad
     && \lm=\pm 2,\nn\\
      N_1(\lm,b)=\left(\begin{array}{cc}\lm&b\\0&\lm\end{array}\right),\qquad
     && \lm=\pm1,\,b=\pm1,\,0,\nn\\
      R(\theta)=\left(\begin{array}{cc}\cos \theta &-\sin \theta \\
      \sin \theta &\cos \theta \end{array} \right), \qquad
      &&\theta\in(0,\pi)\cup(\pi,2\pi),\nn\\
      N_2(\om,b)=\left(\begin{array}{cc}R(\theta)&  b \\  0 &R(\theta)\end{array}\right),
        \qquad
      &&\theta\in(0,\pi)\cup(\pi,2\pi),\nn\eea
where $b= \left(\begin{array}{cc}b_1&b_2\\b_3&b_4\end{array}\right)$
with $b_i\in\R$ and $b_2\neq b_3$.

For any $M\in \Sp(2n)$ and $\om\in\U$, {\it splitting number} of $M$
at $\om$ is defined by
   \bea S_M^{\pm}(\om)=\lim_{\epsilon\to 0^+}
   i_{\om\exp(\pm\sqrt{-1}\epsilon)}(\ga)-i_\om(\ga)\nn\eea
for any path $\ga\in\mathcal{P}_\tau(2n)$ satisfying $\ga(\tau)=M$.

Splitting numbers possesses the following properties.

\noindent{\bf Lemma 2.1.} (cf. \cite{Long0}, Lemma 9.1.5 and List
9.1.12 of \cite{Long1}) {\it Splitting number $S_M^{\pm}(\om)$ are
well defined, i.e., they are independent of the choice of the path
$\ga\in\mathcal{P}_\tau(2n)$ satisfying $\ga(\tau)=M$. For
$\om\in\U$ and $M\in\Sp(2n)$, $S_Q^{\pm}(\om)=S_M^{\pm}(\om)$ if
$Q\approx M$. Moreover we have

(1) $(S_M^+(\pm1),S_M^-(\pm1))=(1,1)$ for $M=\pm N_1(1, b)$ with
$b=1$ or $0$;

(2) $(S_M^+(\pm1),S_M^-(\pm1))=(0,0)$ for $M=\pm N_1(1, b)$ with
$b=-1$;

 (3) $(S_M^+(e^{\sqrt{-1}\theta}),S_M^-(e^{\sqrt{-1}\theta}))=(0,1)$ for
   $M=R(\theta)$ with $\theta\in(0,\pi)\cup (\pi,2\pi)$;

(4) $(S_M^+(\om),S_M^-(\om))=(0,0)$
 for $\om\in \U\setminus \R$ and $M=N_2(\om, b)$ is {\bf trivial} i.e., for

  sufficiently small
 $\alpha>0$, $MR((t-1)\alpha)^{\diamond n}$ possesses no eigenvalues
 on $\U$ for $t\in [0,1)$.

(5) $(S_M^+(\om),S_M^-(\om)=(1,1)$
 for $\om\in \U\setminus \R$ and $M=N_2(\om, b)$ is {\bf non-trivial}.

 (6) $(S_M^+(\om),S_M^-(\om)=(0,0)$
for any $\om\in\U$ and $M\in \Sp(2n)$ with $\sg(M)\cap
\U=\emptyset$.

 (7) $S_{M_1\diamond M_2}^\pm(\om)=S_{M_1}^\pm(\om)+S_{M_2}^\pm(\om)$, for any $M_j\in \Sp(2n_j)$ with
    $j=1,2$ and $\om\in\U$.}

Let \bea F=\R^{2n}\oplus \R^{2n}\nn\eea possess the standard inner
product. We define the symplectic structure of $F$ by
      \bea \{v,w\}=(\mathcal{J}v,w),\;\forall v,w\in F,\;
 {\rm where}\; \mathcal{J}=(-J)\oplus J=\left(\begin{array}{cc} -J &0\\0&J\end{array}\right).
 \;\nn\eea
  We denote by $\Lag(F)$ the set of Lagrangian subspaces of $F$,
  and equip it with the topology as a subspace of the Grassmannian of all
  $ 2n$-dimensional subspaces of $F$.

 It is easy to check that, for any $M\in \Sp(2n)$ its
  graph
      $$\Gr(M)\equiv\left\{\left(\begin{array}{c}x\\Mx\end{array}\right)|x\in
      \R^{2n}\right\}$$
is a Lagrangian subspace of $F$.

Let \bea V_1=\{0\}\times \R^n\times \{0\}\times \R^n\subset \R^{4n},
\quad V_2=\R^n\times \{0\}\times \R^n\times\{0\}\subset
\R^{4n}.\nn\eea

By Proposition 6.1 of \cite{LZ} and Lemma 2.8 and Definition 2.5 of
\cite{LZZ}, we give the following definition.

 \noindent{\bf Definition 2.3.} For any continuous path $\ga\in\mathcal{P}_\tau(2n)$, we
 define the following Maslov-type indices:

 \bea && i_{L_0}(\ga)=\mu^{CLM}_{F}(V_1, \Gr(\ga),[0,\tau])-n,\nn\\
     && i_{L_1}(\ga)=\mu^{CLM}_{F}(V_2, \Gr(\ga),[0,\tau])-n,\nn\\
      && \nu_{L_j}(\ga)=\dim (\ga(\tau)L_j\cap L_j),\qquad j=0,1,\nn\eea
where we denote by $i^{CLM}_F(V,W,[a,b])$ the Maslov index for
Lagrangian subspace path pair $(V,W)$ in $F$ on $[a,b]$ defined by
Cappell, Lee, and Miller in \cite{CLM}. For any $M\in\Sp(2n)$ and
$j=0,1$, we also denote by $\nu_{L_j}(M)=\dim(ML_j\cap L_j)$.

\noindent{\bf Definition 2.4.}  For two paths
$\gamma_0,\;\gamma_1\in
 \mathcal{P_\tau}(2n)$ and $j=0,1$, we say that they are $L_j$-homotopic and denoted by
 $\gamma_0\sim_{L_j}\gamma_1$, if there is a continuous map $\delta:[0,1]\to
 \mathcal{P}(2n)$ such that $\delta(0)=\gamma_0$ and $\delta(1)=\gamma_1$, and
 $\nu_{L_j}(\delta(s))$ is constant for $s\in [0,1]$.

  \noindent{\bf Lemma 2.2.}(\cite{Liu2}) {\it (1) If  $\gamma_0\sim_{L_j}\gamma_1$, there hold
  $$i_{L_j}(\gamma_0)=i_{L_j}(\gamma_1),\;\nu_{L_j}(\gamma_0)=\nu_{L_j}(\gamma_1).$$
 (2) If $\gamma=\gamma_1\diamond \gamma_2\in
  \mathcal{P}(2n)$, and correspondingly $L_j=L_j'\oplus L_j''$, then
  $$i_{L_j}(\gamma)=i_{L'_j}(\gamma_1)+i_{L_j''}(\gamma_2),\;\nu_{L_j}(\gamma)=\nu_{L'_j}(\gamma_1)+\nu_{L_j''}(\gamma_2).$$
 (3) If
$\gamma\in \mathcal{P}(2n)$ is the fundamental solution of
$$\dot x(t)=JB(t)x(t)$$ with symmetric matrix function
$B(t)=\pmatrix b_{11}(t) & b_{12}(t)\\b_{21}(t) & b_{22}(t)
\endpmatrix$ satisfying $b_{22}(t)>0$ for any $t\in R$, then there holds
$$i_{L_0}(\gamma)=\sum_{0<s<1}\nu_{L_0}(\gamma_s),\;\gamma_s(t)=\gamma(st).$$
(4) If $b_{11}(t)>0$ for any $t\in \R$, there holds
$$i_{L_1}(\gamma)=\sum_{0<s<1}\nu_{L_1}(\gamma_s),\;\gamma_s(t)=\gamma(st).$$  }

\noindent {\bf Definition 2.5.} For any $\ga\in
\mathcal{P}_\tau$ and $k\in \N\equiv\{1,2,...\}$, in this paper the
$k$-time iteration $\ga^k$ of $\ga\in \mathcal{P}_\tau(2n)$ in brake
orbit boundary sense is defined by $\td{\ga}|_{[0,k\tau]}$ with \bea
   \td{\ga}(t)=\left\{\begin{array}{l}\ga(t-2j\tau)(N\ga(\tau)^{-1}N\ga(\tau))^j,
   \; t\in[2j\tau,(2j+1)\tau], j=0,1,2,...\\
      N\ga(2j\tau+2\tau-t)N(N\ga(\tau)^{-1}N\ga(\tau))^{j+1},\; t\in[(2j+1)\tau,(2j+2)\tau],
      j=0,1,2,...\end{array}\right.\nn\eea
By \cite{LZZ} or Corollary 5.1 of \cite{LiuZhang} $\displaystyle\lim_{k\to
\infty}\frac{i_{L_0}(\gamma^k)}{k}$ exists, as usual we define the
mean $i_{L_0}$ index of $\ga$ by $\hat{i}_{L_0}(\ga)=\displaystyle\lim_{k\to
\infty}\frac{i_{L_0}(\gamma^k)}{k}$.

For any $P\in\Sp(2n)$ and $\varepsilon\in\R$, we set
 \bea M_\varepsilon(P)=P^T\left(\begin{array}{cc}\sin{2\var}I_n&-\cos{2\var I_n}\\-\cos{2\var}I_n&-\sin
 2\var I_n
 \end{array}\right)P+\left(\begin{array}{cc}\sin{2\var}I_n&\cos{2\var}I_n\\\cos{2\var}I_n&-\sin2\var
 I_n
 \end{array}\right).\nn\eea
Then we have the following

\noindent{\bf Theorem 2.1.}(Theorem 2.3 of \cite{Zhang1})  {\it For
$\ga\in\mathcal{P}_\tau(2k)$ with $\tau>0$, we have
   \bea i_{L_0}(\ga)-i_{L_1}(\ga)=\frac{1}{2}\sgn
M_\var(\ga(\tau)),\nn\eea where $\sgn
M_\var(\ga(\tau))=m^+(M_\var(\ga(\tau)))-m^-(M_\var(\ga(\tau)))$ is
the signature of the symmetric matrix $M_\var(\ga(\tau))$ and
$0<\var\ll 1$. we also have,
 \bea
 (i_{L_0}(\ga)+\nu_{L_0}(\ga))-(i_{L_1}(\ga)+\nu_{L_1}(\ga))=\frac{1}{2}\sign
M_\var(\ga(\tau)),\nn\eea where $0<-\var\ll 1$.}

\noindent{\bf Remark 2.2.} (Remark 2.1 of \cite{Zhang1}) {\it For
any $n_j \times n_j$ symplectic matrix $P_j$ with $j=1,2$ and
 $n_j\in\N$, we have
    \bea &&M_\var(P_1\diamond P_2)=M_\var(P_1)\diamond
    M_\var(P_2),\nn\\
     &&\sgn M_\var(P_1\diamond P_2)=\sgn M_\var(P_1)+
    \sgn M_\var(P_2),\nn\eea
  where $\var\in\R$.}

In the following of this section we will give some lemmas which will
be used frequently in the proof of our main theorem later.

\noindent{\bf Lemma 2.3.} {\it For $k\in\N$ and any symplectic
matrix $P=\left(\begin{array}{cc}I_k&0\\C&I_k\end{array}\right)$,
there holds $P\approx I_2^{\diamond p}\diamond N_1(1,1)^{\diamond
q}\diamond N_1(1,-1)^{\diamond r}$ with $p,q,r$  satisfying
  \bea m^0(C)=p, \quad m^-(C)=q,\quad m^+(C)=r.\nn\eea
  }

 {\bf Proof.} It is clear that
      \bea P\approx
      \left(\begin{array}{cc}I_k&0\\B&I_k\end{array}\right),\nn\eea
  where $B=\diag (0,-I_{m^-(C)},I_{m^+(C)})$. Since
$J_1N_1(1,\pm
1)(J_1)^{-1}=\left(\begin{array}{cc}1&0\\\mp1&1\end{array}\right)$,
by Remark 2.1 we have $N_1(1,\pm 1) \approx
 \left(\begin{array}{cc}1&0\\\mp1&1\end{array}\right)$. Then
  \bea P\approx I_2^{\diamond m^0(C)}\diamond N_1(1,1)^{\diamond m^-(C)}\diamond N_1(1,-1)^{\diamond
m^+(C)}.\nn\eea

By Lemma 2.1 we have
  \be S_P^+(1)=m^0(C)+m^-(C)=p+q.\lb{z1}\ee
By the definition of the relation $\approx$, we have
 \be 2p+q+r=\nu_1(P)=2m^0(C)+m^+(C)+m^-(C).\lb{z_2}\ee
Also we have
 \be p+q+r=m^0(C)+m^+(C)+m^-(C)=k.\lb{zh3}\ee
By (\ref{z1})-(\ref{zh3}) we have
 \bea m^0(C)=p, \quad m^-(C)=q,\quad m^+(C)=r.\nn\eea
The proof of Lemma 2.3 is complete.\hfill\hb

\noindent{\bf Definition 2.6.} We call two symplectic matrices $M_1$
and $M_2$ in $\Sp(2k)$ are {\it special homotopic}(or
$(L_0,L_1)$-homotopic) and denote by $M_1\sim M_2$, if there are
$P_j\in\Sp(2k)$ with $P_j=\diag (Q_j,(Q_j^T)^{-1})$, where $Q_j$ is
a $k\times k$ invertible real matrix, and $\det(Q_j)>0$ for $j=1,2$,
such that
    $$M_1=P_1M_2P_2.$$
It is clear that $\sim$ is an equivalent relation.

\noindent{\bf Lemma 2.4.}  {\it For $M_1,\,M_2\in \Sp(2k)$, if
$M_1\sim M_2$, then
 \bea
        {\r \sgn M_\var(M_1)=\sgn M_\var(M_2),\quad 0\le |\var|\ll 1,\lb{z6}}\\
       N_kM_1^{-1}N_kM_1\approx N_kM_2^{-1}N_kM_2.\lb{z7}\eea}

 {\bf Proof.}
By Definition 2.6, there are $P_j\in\Sp(2k)$ with $P_j=\diag
(Q_j,(Q_j^T)^{-1})$, $Q_j$ being $k\times k$ invertible real matrix,
and $\det(Q_j)>0$ such that
    $$M_1=P_1M_2P_2.$$
 Since $\det(Q_j)>0$ for $j=1,2$, we can joint $Q_j$ to
  $I_k$ by invertible matrix path. Hence we can joint $P_1M_2P_2$ to
  $M_2$ by symplectic path preserving the nullity $\nu_{L_0}$ and $\nu_{L_1}$.
  By Lemma 2.2 of \cite{Zhang1}, (\ref{z6})
  holds. Since $P_jN_k=N_kP_j$ for $j=1,2$. Direct computation shows that
    \be N_k(P_1M_2P_2)^{-1}N_k(P_1M_2P_2)=P_2^{-1}N_kM_2^{-1}N_kM_2P_2.\ee
   Thus (\ref{z7}) holds from Remark 2.1. The proof of Lemma 2.4 is
   complete.\hfill\hb

\noindent{\bf Lemma 2.5.} {\it Let
$P=\left(\begin{array}{cc}A&B\\C&D\end{array}\right)\in \Sp(2k)$,
where $A,B,C,D$ are all $k\times k$ matrices. Then

(i) $\frac{1}{2}\sgn M_\var(P)\le k-\nu_{L_0}(P)$, for $0<\var\ll
1$. If $B=0$, we have $\frac{1}{2}\sgn M_\var(P)\le 0$ for
$0<\var\ll 1$.

(ii) Let $m^+(A^TC)=q$, we have
     \bea \frac{1}{2}\sgn M_\var(P)\le k-q,\quad 0\le|\var|\ll 1.\lb{zhang5}\eea
     Moreover if $B=0$, we have
      \bea \frac{1}{2}\sgn M_\var(P)\le -q,\quad 0<-\var\ll 1.\lb{zhang6}\eea

(iii) $\frac{1}{2}\sgn M_\var(P)\ge \dim\ker C-k$ for
      $0<\var\ll 1$, If $C=0$, then $\frac{1}{2}\sgn M_\var(P)\ge 0$ for $0<\var\ll 1$

(iv)  If both $B$ and $C$ are invertible, we have
       \bea \sgn M_\var(P)=\sgn M_0(P),\quad 0\le |\var|\ll 1.\nn\eea}

{\bf Proof.}
  Since $P$ is symplectic,  so is for $P^T$. From $P^TJ_kP=J_k$ and $PJ_kP^T=J_k$ we get
   $A^TC, B^TD,AB^T,CD^T$ are all symmetric matrices and
   \be AD^T-BC^T=I_k,\quad A^TD-C^TB=I_k. \lb{zhang3}\ee
 We denote by $s=\sin 2\var$ and $c=\cos 2\var$. By
  definition of $M_\var(P)$, we have
   \bea
   M_\var(P)&=&\left(\begin{array}{cc}A^T&C^T\\B^T&D^T\end{array}\right)\left(\begin{array}{cc}sI_k&-cI_k\\-cI_k&-sI_k\end{array}\right)
   \left(\begin{array}{cc}A&B\\C&D\end{array}\right)+\left(\begin{array}{cc}sI_k&cI_k\\cI_k&-sI_k\end{array}\right)\nn\\
          &=& \left(\begin{array}{cc}A^T&C^T\\B^T&D^T\end{array}\right)\left(\begin{array}{cc}sI_k&-2cI_k\\0&-sI_k\end{array}\right)
   \left(\begin{array}{cc}A&B\\C&D\end{array}\right)+\left(\begin{array}{cc}sI_k&2cI_k\\0&-sI_k\end{array}\right)\nn\\
         &=& \left(\begin{array}{cc}sA^TA-2cA^TC-sC^TC+sI_k& *
         \\sB^TA-2cB^TC-sD^TC&sB^TB-2cB^TD-sD^TD-sI_k\end{array}\right)\nn\\
         &=& \left(\begin{array}{cc}sA^TA-2cA^TC-sC^TC+sI_k&
         sA^TB-2cC^TB-sC^TD
         \\sB^TA-2cB^TC-sD^TC&sB^TB-2cB^TD-sD^TD-sI_k\end{array}\right),\lb{zhang1}\eea
where in the second equality we have used that $P^TJ_kP=J_k$, in the
fourth equality we have used that $M_\var(P)$ is a symmetric matrix.
    So
    \bea M_0(P)=-2\left(\begin{array}{cc}A^TC&C^TB\\B^TC&B^TD\end{array}\right)
    = -2\left(\begin{array}{cc}C^T&0\\0&B^T\end{array}\right)\left(\begin{array}{cc}A&B\\C&D\end{array}\right),\nn\eea
     where we have used $A^TC$ is symmetric. So if both $B$ and $C$ are invertible, $M_0(P)$ is invertible and symmetric, its signature is invariant under
      small perturbation, so (iv) holds.

       If $\nu_{L_0}(P)=\dim\ker B>0$, since $B^TD=D^TB$, for any $x\in \ker
       B\subseteq\R^k$, $x\neq 0$, and $0<\var\ll 1$, we have
        \bea
        &&M_\var(P)\left(\begin{array}{c}0\\x\end{array}\right)\cdot \left(\begin{array}{c}0\\x\end{array}\right)
        =(sB^TB-2cD^TB-sD^TD-sI_k)x\cdot x\nn\\
         &&=-s(D^TD+I_k)x\cdot x\nn\\
         &&<0. \lb{zhang4}\eea
       So $M_\var(P)$ is negative definite on $(0\oplus \ker
       B)\subseteq \R^{2k}$. Hence $m^-(M_\var(p)\ge \dim\ker B$
       which yields that $\frac{1}{2}\sgn M_\var(P)\le k-\dim\ker
       B=k-\nu_{L_0}(P)$, for $0<\var\ll 1$. Thus (i) holds. Similarly we can prove (iii).

       If $m^+(A^TC)=q>0$, let $A^TC$ is positive definite on
       $E\subseteq\R^k$, then for $0\le |s|\ll 1$, similar to
       (\ref{zhang4}) we have $M_\var(P)$ is negative on $E\oplus
       0\subseteq \R^{2k}$. Hence $m^-(M_\var(P)\ge q$, which yields
       (\ref{zhang5}).

       If $B=0$, by (\ref{zhang1}) we have
       \bea M_\var(P)= \left(\begin{array}{cc}sA^TA-2cA^TC-sC^TC+sI_k&
        -sC^TD
         \\-sD^TC&-sD^TD-sI_k\end{array}\right).\eea

    Since
        \bea && \left(\begin{array}{cc}I_k&-C^TD(D^TD+I_k)^{-1}\\0&I_k\end{array}\right)
              \left(\begin{array}{cc}sA^TA-2cA^TC-sC^TC+sI_k&-sC^TD
              \\-sD^TC&-sD^TD-sI_k\end{array}\right)\cdot\nn\\
              &&\quad\cdot\left(\begin{array}{cc}I_k&0\\-(D^TD+I_k)^{-1}D^TC&I_k\end{array}\right)\nn\\
              &&=
              \left(\begin{array}{cc}sA^TA-2cA^TC-sC^TC+sI_k+sC^TD(D^TD+I_k)^{-1}D^TC
              &0\\0&-sD^TD-sI_k\end{array}\right),\eea
     for $0< -s\ll 1$, we have
     \be m^-(M_\var(P))\ge k+m^+(A^TC)\ee
     which yields (\ref{zhang6}). So (ii) holds and the proof of
     Lemma 2.5 is complete.\hfill\hb

\noindent{\bf Lemma 2.6.} (\cite{Zhang1}) {\it For
$\ga\in\mathcal{P}_\tau(2)$, $b>0$, and $0<\var\ll 1$ small enough
we have
 \bea
&& \sgn M_{\pm\var}(R(\theta))=0,\quad
 {\rm for}\;\theta\in \R,\nn\\
&&  \sgn M_\var(P)=0,\quad {\rm if}\; P=\pm
 \left(\begin{array}{cc}1&b\\0&1\end{array}\right)\;{\rm or}\;
 \pm
 \left(\begin{array}{cc}1&0\\-b&1\end{array}\right),\nn\\
&&\sgn M_\var(P)=2,\quad {\rm if}\; P=\pm
 \left(\begin{array}{cc}1&-b\\0&1\end{array}\right),\nn\\
 &&\sgn M_\var(P)=-2,\quad {\rm if}\; P=
 \pm
 \left(\begin{array}{cc}1&0\\b&1\end{array}\right).\nn\eea}

\setcounter{equation}{0}
\section{Proofs of Theorems 1.1 and 1.2.}
In this section we prove Theorems 1.1 and 1.2. The proof mainly
depends on the method in \cite{LiuZhang} and the following

\noindent{\bf Theorem 3.1.} {\it For any odd number $n\ge 3$,
$\tau>0$ and $\ga\in \P_\tau(2n)$, let $P=\ga(\tau)$. If $i_{L_0}\ge
0$, $i_{L_1}\ge 0$, $i(\ga)\ge n$, $\ga^2(t)=\ga(t-\tau)\ga(\tau)$
for all $t\in[\tau,2\tau]$, and $P\sim (-I_2)\diamond Q$ with
$Q\in\Sp(2n-2)$, then
 \be i_{L_1}(\ga)+S_{P^2}^+(1)-\nu_{L_0}(\ga)>\frac{1-n}{2}.\lb{dd}\ee}

{\bf Proof.} If the conclusion of Theorem 3.1 does not hold, then
        \be i_{L_1}(\ga)+S_{P^2}^+(1)-\nu_{L_0}(\ga)\le
        \frac{1-n}{2}.\lb{d0}\ee
  In the following we shall obtain a contradiction from (\ref{d0}). Hence (\ref{dd})
  holds and Theorem 3.1 is proved.

Since $n\ge 3$ and $n$ is odd, in the following of the proof of
Theorem 3.1 we write $n=2p+1$ for some $p\in\N$. We denote by
$Q=\left(\begin{array}{cc} A&B\\C&D\end{array}\right)$, where
$A,B,C,D$ are $(n-1)\times (n-1)$ matrices. Then since $Q$ is a
symplectic matrix we have
    \be A^TC=C^TA,\;B^TD=D^TB,\; AB^T=BA^T,\;CD^T=DC^T,\lb{zz1}\ee
     \be AD^T-BC^T=I_{n-1},\quad A^TD-C^TB=I_{n-1},\lb{zz2}\ee
     \be \dim\ker B=\nu_{L_0}(\ga)-1,\quad \dim\ker
     C=\nu_{L_1}(\ga)-1.\lb{zz3}\ee
 Since $\ga^2(t)=\ga(t-\tau)\ga(\tau)$ for all $t\in[\tau,2\tau]$
we have $\ga^2$ is also the twice iteration of $\ga$ in the periodic
boundary value case, so by the Bott-type formula (cf. Theorem 9.2.1
of \cite{Long1}) and the proof of Lemma 4.1 of \cite{LZZ} we have
 \bea   && i(\ga^2)+2S_{P^2}^+(1)-\nu(\ga^2)\nn\\
       &=& 2i(\ga)+2S_{P}^+(1)+\sum_{\theta\in
       (0,\pi)}(S_{P}^+(e^{\sqrt{-1}\theta})\nn\\
       &&-(\sum_{\theta\in
       (0,\pi)}(S_{P}^-(e^{\sqrt{-1}\theta})+(\nu(P)-S_{P}^-(1))+(\nu_{-1}(P)-S_{P}^-(-1)))\nn\\
       &\ge& 2n+2S_{P}^+(1)-n\nn\\
       &=&n+2S_{P}^+(1)\nn\\
        &\ge& n, \lb{d1} \eea
where we have used the condition $i(\ga)\ge n$ and
$S^+_{P^2}(1)=S^+_P(1)+S^+_P(-1)$,
$\nu(\gamma^2)=\nu(\gamma)+\nu_{-1}(\gamma)$. By by Proposition C of
\cite{LZZ} and Proposition 6.1 of \cite{LiuZhang} we have
 \be i_{L_0}(\ga)+i_{L_1}(\ga)=i(\ga^2)-n,\quad \nu_{L_0}(\ga)+\nu_{L_1}(\ga)=\nu(\ga^2).\lb{d2}\ee
So by (\ref{d1}) and (\ref{d2}) we have
 \bea
 &&(i_{L_1}(\ga)+S_{P^2}^+(1)-\nu_{L_0}(\ga))+(i_{L_0}(\ga)+S_{P^2}^+(1)-\nu_{L_1}(\ga))\nn\\
 &&=i(\ga^2)+2S_{P^2}^+(1)-\nu(\ga^2)-n\nn\\
 &&\ge n-n\nn\\
 &&=0.\lb{d3}\eea
By Theorem 2.1 and Lemma 2.6 we have
 \bea
 &&(i_{L_1}(\ga)+S_{P^2}^+(1)-\nu_{L_0}(\ga))-(i_{L_0}(\ga)+S_{P^2}^+(1)-\nu_{L_1}(\ga))\nn\\
 &&=i_{L_1}(\ga)-i_{L_0}(\ga)-\nu_{L_0}(\ga))+\nu_{L_1}(\ga)\nn\\
 &&=-\frac{1}{2}\sgn M_\var(Q)-\frac{1}{2}\sgn M_\var(-I_2)\nn\\
 &&=-\frac{1}{2}\sgn M_\var(Q)\nn\\
 &&\ge 1-n.\lb{d4}\eea
So by (\ref{d3}) and (\ref{d4}) we have
 \be i_{L_1}(\ga)+S_{P^2}^+(1)-\nu_{L_0}(\ga)\ge
 \frac{1-n}{2}.\lb{d5}\ee

By (\ref{d0}), the inequality of (\ref{d5}) must be equality. Then
both (\ref{d1}) and (\ref{d4}) are equality. So we have

 \be i(\ga^2)+2S_{P^2}^+(1)-\nu(\ga^2)=n.\lb{d7}\ee
 \be  i_{L_1}(\ga)+S_{P^2}^+(1)-\nu_{L_0}(\ga)=
 \frac{1-n}{2}.\lb{d8}\ee
 \be
 i_{L_0}(\ga)+\nu_{L_0}(\ga)-i_{L_1}(\ga)-\nu_{L_1}(\ga)=n-1.\lb{d9}\ee

Thus by (\ref{d1}), (\ref{d7}),  Theorem 1.8.10 of \cite{Long1}, and
Lemma 2.1 we have
  \bea P\approx (-I_2)^{\diamond p_1}\diamond N_1(1,-1)^{\diamond p_2}
       \diamond N_1(-1,1)^{\diamond p_3}\diamond
       R(\theta_1)\diamond R(\theta_2)\diamond\cdots\diamond
       R(\theta_{p_4}),\nn \eea
where $p_j\ge 0$ for $j=1,2,3,4$, $p_1+p_2+p_3+p_4=n$ and
$\theta_j\in(0,\pi)$ for $1\le j\le p_4$. Otherwise by
(\ref{d1}) and Lemma 2.1 we have
$i(\ga^2)+2S_{P^2}^+(1)-\nu(\ga^2)>n$ which contradicts to
(\ref{d7}). So by Remark 2.1, we have \be P^2 \approx I_2^{\diamond
p_1}\diamond N_1(1,-1)^{\diamond
      p_2}\diamond
       R(\theta_1)\diamond R(\theta_2)\diamond\cdots\diamond
       R(\theta_{p_3}),\lb{d10} \ee
where $p_i\ge 0$ for $1\le i\le 3$, $p_1+p_2+p_3=n$ and $\theta_j\in
(0,2\pi)$ for $1\le j\le p_3$.

      Note that, since $\ga^2(t)=\ga(t-\tau)\ga(\tau)$, we have
           \bea \ga^2(2\tau)=\ga(\tau)^2=P^2.\lb{g1}\eea
       By Definition 2.5 we have
            \bea \ga^2(2\tau)=N\ga(\tau)^{-1}N\ga(\tau)=NP^{-1}NP.\lb{g2}\eea
   So by (\ref{g1}) and (\ref{g2}) we have
          \bea P^2=NP^{-1}NP.\lb{g3}\eea
        By (\ref{g3}), Lemma 2.4, and $P\sim (-I_2)\diamond Q$ we have
             \bea P^2&=&NP^{-1}NP\nn\\
                     &\approx& N((-I_2)\diamond Q)^{-1}N((-I_2)\diamond
                     Q)\nn\\
                     &=&I_2\diamond (N_{n-1}Q^{-1}N_{n-1}Q).\lb{g4}\eea
           So by (\ref{d10}), we have
                \bea p_1\ge 1.\lb{g5}\eea
Also by (\ref{g4}) and Lemma 2.5, we have
   \bea P^2\approx I_2\diamond (N_{n-1}Q'^{-1}N_{n-1}Q'),\quad
   \forall\, Q'\sim Q\; {\rm where}\; Q'\in\Sp(2n-2).\lb{g6}\eea

By (\ref{d10}) it is easy to check that

\be {\rm tr} (P^2)=2n-2p_3+2\sum_{j=1}^{p_3}\cos
\theta_j.\lb{d20}\ee

By (\ref{d7}), (\ref{d10}) and Lemma 2.1 we have
 \bea n=i(\ga^2)+2S_{P^2}^+(1)-\nu(\ga^2)=i(\ga^2)-p_2\ge
 i(\ga^2)-n+1.\nn\eea
So
 \be i(\ga^2)\le 2n-1.\lb{11}\ee
 By (\ref{d2}) we have
  \be i(\ga^2)=n+i_{L_0}(\ga)+i_{L_1}(\ga).\lb{d12}\ee
Since $i_{L_0}(\ga)\ge0$ and $i_{L_1}(\ga)\ge 0$, we have $n\le
i(\ga^2)\le 2n-1$. So we can divide the index $i(\ga^2)$ into the
following three cases.

{\bf Case I.}  $i(\ga^2)=n$.

In this case, by (\ref{d2}), $i_{L_0}(\ga)\ge 0$, and
$i_{L_1}(\ga)\ge 0$, we have
 \be i_{L_0}(\ga)=0=i_{L_1}(\ga).\lb{d13}\ee
So by (\ref{d9}) we have
 \be \nu_{L_0}(\ga)-\nu_{L_1}(\ga)=n-1.\ee
 Since $\nu_{L_1}(\ga)\ge 1$ and $\nu_{L_0}(\ga)\le n$, we have
  \be \nu_{L_0}(\ga)=n,\quad \nu_{L_1}(\ga)=1.\lb{d14}\ee
By (\ref{d2}) we have
 \be\nu(\ga^2)=\nu(P^2)=n+1.\lb{d60}\ee
 By (\ref{d8}), (\ref{d13}) and (\ref{d14}) we have
 \be S_{P^2}^+(1)=\frac{1-n}{2}+n=\frac{1+n}{2}=p+1.\lb{d61}\ee
 So by (\ref{d10}), (\ref{d60}), (\ref{d61}), and Lemma 2.1 we have
  \be P^2\approx I_2^{\diamond (p+1)}\diamond
  R(\theta_1)\diamond\cdots\diamond R(\theta_p),\lb{d90}\ee
  where $\theta_j\in (0,2\pi)$. By (\ref{zz3})  and (\ref{d14}) we have
  $B=0$. By (\ref{g4}), (\ref{zz1}), and (\ref{zz2}), we have
   \bea P^2&=&NP^{-1}NP\approx I_2\diamond (N_{n-1}Q^{-1}N_{n-1}Q)\nn\\
   &=&I_2\diamond\left(\begin{array}{cc} D^T&0\\C^T&A^T\end{array}\right)
   \left(\begin{array}{cc} A&0\\C&D\end{array}\right)\nn\\&=&
   I_2\diamond\left(\begin{array}{cc} D^TA&0\\2C^TA&AD^T\end{array}\right)\nn\\
   &=&I_2\diamond \left(\begin{array}{cc}
   I_{2p}&0\\2A^TC&I_{2p}\end{array}\right).\nn\eea
   Hence $\sg(P^2)=\{1\}$ which contradicts to (\ref{d90}) since $p\ge 1$.

{\bf Case II.} $i(\ga^2)=n+2k$, where $1\le k\le p$.

  In this case by (\ref{d2}) we have
   \bea i_{L_0}(\ga)+i_{L_1}(\ga)=2k.\nn\eea
   Since $i_{L_0}(\ga)\ge 0$ and $i_{L_1}(\ga)\ge 0$  we can write $i_{L_0}(\ga)=k+r$ and $i_{L_1}(\ga)=k-r$ for
   some integer $-k\le r \le k$. Then by (\ref{d9}) we have
    \be n-1\ge \nu_{L_0}(\ga)-\nu_{L_1}(\ga)=n-2r-1.\lb{d15}\ee
    Thus $r\ge 0$ and $0\le r\le k$.

  By Theorem 2.1 and (i) of Lemma 2.5 we have
   \be 2r=i_{L_0}(\ga)-i_{L_1}(\ga)=\frac{1}{2}M_\var (P)\le
   n-\nu_{L_0}(P)\ee
   which yields that $\nu_{L_0}(\ga)\le n-2r$. So by (\ref{d15}) and $\nu_{L_1}(\ga)\ge 1$ we have
      \be \nu_{L_0}(\ga)=n-2r,\quad \nu_{L_1}(\ga)=1.\lb{d35}\ee
   Then by (\ref{d8}) we have
    \be
    S_{P^2}^+(1)=(n-2r)+\frac{1-n}{2}-(k-r)=\frac{1+n}{2}-k-r=p+1-k-r.\lb{d23}\ee
  Then by (\ref{d10}) and $\nu(P^2)=n-2r+1$ and Lemma 2.1 we have
  \be P^2\approx I_2^{\diamond (p+1-k-r)}\diamond
  N_1(1,-1)^{\diamond 2k}\diamond
     R(\theta_1)\diamond \cdots\diamond R(\theta_q),\ee
     where $q=n-(p+1-k-r)-2k=p+r-k\ge 0$.
Then we have the following three subcases (i)-(iii).

 {\bf (i) $q=0$}.

    The only possibility is $k=p$ and $r=0$, in this case
    $P^2\approx I_2\diamond N_1(1,-1)^{\diamond 2p}$ and $B=0$.
     By direct computation we have
     \bea N_1(1,-1)^{\diamond 2p}\approx N_{2p}Q^{-1}N_{2p}Q=\left(\begin{array}{cc}
   I_{n-1}&0\\2A^TC&I_{n-1} \end{array}\right).\lb{d65}\eea
   Then by Lemma 2.3 we have
   \bea m^+(A^TC)=2p.\nn\eea
   By (ii) of Lemma 2.5 we have
   \be \frac{1}{2}\sgn M_\var(Q)\le 2p-2p=0,\qquad 0< -\var\ll 1.\lb{d17}\ee
   Thus by (\ref{d17}) and Theorem 2.1, for  $0< -\var\ll 1$ we have,
   \bea
   &&(i_{L_0}(\ga)+\nu_{L_0}(\ga))-(i_{L_1}(\ga)+\nu_{L_1}(\ga))\nn\\
   &=&\frac{1}{2}\sgn M_\var(P)\nn\\
   &=&\frac{1}{2}\sgn M_\var (I_2)+\frac{1}{2}M_\var(Q)\nn\\
   &=&0+\frac{1}{2}M_\var(Q)\nn\\
   &\le& 0\nn\eea
  which contradicts (\ref{d9}).

{\bf (ii) $q>0$ and $r=0$}.

 In this case $\nu_{L_0}(\ga)=n$ and $\nu_{L_1}(\ga)=1$, also we have $B=0$.
By the equality of (\ref{d65}) we have
   \bea {\rm tr}\,(P^2)=2n \nn\eea
  which contradicts to (\ref{d20}) with $p_3=q>0$.

{\bf (iii) $q>0$ and $r>0$}.

In this case, by (\ref{d23}) we have $r<p$ (otherwise, then $p=r=k$.
From (\ref{g5}) there holds $S^+_{P^2}(1)\ge 1$, so from (\ref{d23})
we have $1\le S^+_{P^2}(1)=1-p\le 0$ a contradiction). Here it is
easy to see $\rank B=2r$. Then there are two invertible $2p\times
2p$ matrices $U$ and $V$  with $\det U>0$ and $\det V>0$ such that
 \bea UBV= \left(\begin{array}{cc}
   I_{2r}&0\\0&0\end{array}\right).\nn\eea
   So there holds
   \be  Q\sim\, \diag (U,(U^T)^{-1})Q\diag ((V^T)^{-1}, V)=\left(\begin{array}{cccc}
   A_1&B_1&I_{2r}&0\\C_1&D_1&0&0\\A_3&B_3&A_2&B_2\\C_3&D_3&C_2&D_2\end{array}\right):=Q_1,\ee
where for $j=1,2,3$, $A_j$ is a  $2r\times 2r$ matrix, $D_j$ is a
$(2p-2r)\times (2p-2r)$ matrix for $j=1,2,3$, $B_j$ is a $2r\times
(2p-2r)$ matrix, and $C_j$ is $(2p-2r)\times 2r$ matrix. Since $Q_1$
is still a symplectic matrix, we have $Q_1^TJ_{2p}Q_1=J_{2p}$, then
it is easy to check that
 \be C_1=0, \; B_2=0.\ee
So
 \be Q_1=\left(\begin{array}{cccc}
   A_1&B_1&I_{2r}&0\\0&D_1&0&0\\A_3&B_3&A_2&0\\C_3&D_3&C_2&D_2\end{array}\right).\ee
So for the case (iii) of Case II, we have the following 3 subcases
1-3.

{\bf Subcase 1.} $A_3=0$.

     In this case since $Q_1$ is symplectic, by direct computation
     we  have
     \bea N_{2p}Q_1^{-1}N_{2p}Q_1=\left(\begin{array}{cccc}
   I_{2r}&*&*&*\\**&I_{2p-2r}&*&*\\**&*&I_{2r}&*\\**&*&*&I_{2p-2r}\end{array}\right).\nn\eea
   Hence we have
   \bea {\rm tr}(N_{2p}Q_1^{-1}N_{2p}Q_1)=4p.\nn\eea
 Since $Q_1\sim Q$, we have
    \be P\sim (-I_2)\diamond Q_1.\lb{g7}\ee
   Then by the proof of Lemma 2.4 we have
   \bea {\rm tr}P^2&=&{\rm tr} (NP^{-1}NP)\nn\\
        &=&{\rm tr}N((-I_2)\diamond Q_1)^{-1}N((-I_2)\diamond Q_1)\nn\\
                             &=&{\rm tr}\,I_2\diamond ((N_{2p}Q_1^{-1}N_{2p}Q_1)\nn\\
                             &=&4p+2=2n.\lb{d21}\eea

   By (\ref{d20}) and $p_3=q>0$ we have
   \be {\rm tr} (P^2)<2n.\lb{d22}\ee
   (\ref{d21}) and (\ref{d22}) yield a contradiction.

{\bf Subcase 2.} $A_3$ is invertible.

By $Q_1$ is symplectic we have
    \be \left(\begin{array}{cc}
   A^T_1& 0\\B_1^T&D_1^T\end{array}\right)\left(\begin{array}{cc}
   A_2& 0\\C_2&D_2\end{array}\right)-\left(\begin{array}{cc}
   A_3^T&C_3^T\\B_3^T&D_3^T\end{array}\right)\left(\begin{array}{cc}
   I_{2r}&0\\0&0\end{array}\right)=I_{2p}.\ee
Hence
     \be D_1^TD_2=I_{2p-2r}.\lb{d25}\ee

       By direct computation we have
       \bea \left(\begin{array}{cccc}
   A_1&B_1&I_{2r}&0\\0&D_1&0&0\\A_3&B_3&A_2&0\\C_3&D_3&C_2&D_2\end{array}\right)\left(\begin{array}{cccc}
   I_{2r}&-A_3^{-1}B_3&0&0\\0&I_{2p-2r}&0&0\\0&0&I_{2r}&0\\0&0&B_3^T(A_3^T)^{-1}&I_{2p-2r}\end{array}\right)
    =\left(\begin{array}{cccc}
   A_1&\td{B_1}&I_{2r}&0\\0&D_1&0&0\\A_3&0&A_2&0\\C_3&\td{D}_3&\td{C}_2&D_2\end{array}\right).\nn\eea

So by (\ref{d25}) we have
 \bea &&\left(\begin{array}{cccc}
   I_{2r}&-\td{B}_1D_2^T&0&0\\0&I_{2p-2r}&0&0\\0&0&I_{2r} &0\\0&0&D_2\td {B}_1^T&I_{2p-2r}\end{array}\right)\left(\begin{array}{cccc}
   A_1&\td{B_1}&I_{2r}&0\\0&D_1&0&0\\A_3&0&A_2&0\\C_3&\td{D}_3&\td{C}_2&D_2\end{array}\right)\nn\\
  &=&\left(\begin{array}{cccc}
   A_1&0&I_{2r}&0\\0&D_1&0&0\\A_3&0&A_2&0\\\td{C}_3&\td{D}_3&\hat{C}_2&D_2\end{array}\right):=Q_2.\nn \eea
Then we have
 \be Q_2\sim Q_1\sim Q.\lb{d27}\ee

  Since $Q_2$ is a symplectic
matrix, we have $Q_2^TJ_{2p}Q_2=J_{2p}$, then it is easy to check
that
 \be \td{C}_3=0,\; \hat{C}_2=0.\ee
 Hence we have
 \be Q_2= \left(\begin{array}{cc}
   A_1&I_{2r}\\A_3&A_2\end{array}\right)\diamond \left(\begin{array}{cc}
   D_1&0\\\td{D}_3&D_2\end{array}\right).\ee
Since
 \be N_{2p-2r}\left(\begin{array}{cc}
   D_1&0\\\td{D}_3&D_2\end{array}\right)^{-1}N_{2p-2r}\left(\begin{array}{cc}
   D_1&0\\\td{D}_3&D_2\end{array}\right)=\left(\begin{array}{cc}
   I_{2p-2r}&0\\2D_1^T\td{D}_3&I_{2p-2r}\end{array}\right),\ee
by (\ref {d27}), (\ref{g6}),  and Lemma 2.4, there is a symplectic
matrix $W$ such that
   \be P^2\approx I_2\diamond W \diamond \left(\begin{array}{cc}
   I_{2p-2r}&0\\2D_1^T\td{D}_3&I_{2p-2r}\end{array}\right).\ee
Then by (\ref{d10}) and Lemma 2.3, $D_1^T\td{D}_3$ is semipositive
and
 \bea 1+m^0(D_1^T\td{D}_3)\le S_{P^2}^+(1).\nn\eea
 So by (\ref{d23}) we have
  \be m^0(D_1^T\td{D}_3) \le
  p+1-k-r-1=p-k-r=(2p-2r)-(p+k-r)\le 2p-2r-1.\lb{d68}\ee
Since $D_1^T\td{D}_3$ is a semipositive $(2p-2r)\times (2p-2r)$
matrix, by (\ref{d68}) we have $m^+(D_1^T\td{D}_3)>0$. Then by
Theorem 2.1, (ii) of Lemma 2.5 and Lemma 2.6, for $0< -\var\ll 1$ we
have
 \bea
 &&(i_{L_0}(\ga)+\nu_{L_0}(\ga))-(i_{L_1}(\ga)+\nu_{L_1}(\ga))\nn\\
 &=& \frac{1}{2} \left(M_\var(-I_2)+M_\var\left(\left(\begin{array}{cc}
   A_1&I_{2r}\\A_3&A_2\end{array}\right)\right)+M_\var\left(\left(\begin{array}{cc}
   D_1&0\\\td{D}_3&D_2\end{array}\right)\right)\right)\nn\\
 &\le&\frac{1}{2}(0+4r+2(2p-2r-1))\nn\\
  &=&2p-1\nn\\
  &=&n-2\eea
which contradicts to (\ref{d9}).

{\bf Subcase 3.} $A_3\neq 0$ and $A_3$ is not invertible.

In this case, suppose $\rank A_3=\lm$, then $0<\lm<2r$. There is a
invertible $2r\times 2r$ matrix $G$ with $\det G>0$ such that
 \be
 GA_3G^{-1}=\left(\begin{array}{cc}\Lm&0\\0&0\end{array}\right),\lb{d30}\ee
where $\Lm$ is a $\lm\times \lm$ invertible matrix. Then we have
  \bea &&\left(\begin{array}{cccc}
   (G^T)^{-1}&0&0&0\\0&I_{2p-2r}&0&0\\0&0&G&0\\0&0&0&I_{2p-2r}\end{array}\right)\left(\begin{array}{cccc}
   A_1&B_1&I_{2r}&0\\0&D_1&0&0\\A_3&B_3&A_2&0\\C_3&D_3&C_2&D_2\end{array}\right)
   \left(\begin{array}{cccc}
   (G)^{-1}&0&0&0\\0&I_{2p-2r}&0&0\\0&0&G^T&0\\0&0&0&I_{2p-2r}\end{array}\right)\nn\\
    &=& \left(\begin{array}{cccc}
   \td{A_1}&\td{B}_1&I_{2r}&0\\0&D_1&0&0\\GA_3G^{-1}&\td{B}_3&\td{A}_2&0\\  \tilde{C}_3&D_3&  \tilde{C}_2&D_2\end{array}\right):=Q_3.\eea
By (\ref{d30}) we can write $Q_3$ as the following block form
 \bea Q_3=\left(\begin{array}{cccccc}
   U_1&U_2&F_1&I_\lm&0&0\\U_3&U_4&F_2&0&I_{2r-\lm}&0\\0&0&D_1&0&0&0\\\Lm&0&E_1&W_1&W_2&0
   \\0&0&E_2&W_3&W_4&0\\G_1&G_2&D_3&K_1&K_2&D_2\end{array}\right).\lb{d31}\eea
Let $R_1=\left(\begin{array}{ccc}
   I_\lm&0&0\\0&I_{2r-\lm}&0\\-G_1\Lm^{-1}&0&I_{2p-2r}\end{array}\right)$
   and $R_2=\left(\begin{array}{ccc}
   I_\lm&0&-\Lm^{-1}E_1\\0&I_{2r-\lm}&0\\0&0&I_{2p-2r}\end{array}\right)$.
By (\ref{d31}) we have
 \bea \diag ((R_1^T)^{-1},R_1)Q_3\diag(R_2,(R^T_2)^{-1})=
 \left(\begin{array}{cccccc}
   U_1&U_2&\td{F}_1&I_\lm&0&0\\U_3&U_4&\td{F}_2&0&I_{2r-\lm}&0\\0&0&D_1&0&0&0\\\Lm&0&0&W_1&W_2&0
   \\0&0&E_2&W_3&W_4&0\\0&G_2&\td{D}_3&\td{K}_1&\td{K}_2&D_2\end{array}\right):=Q_4.\nn\eea
Since $Q_4$ is a symplectic matrix we have
 \bea Q_4^TJQ_4=J\lb{d33}.\nn\eea
Then by (\ref{d33}) and direct computation we have $U_2=0$, $U_3=0$,
$W_2=0$, $W_3=0$, $\td{F}_1=0$, $\td{K_1}=0$, and $U_1$, $U_4$,
$W_1$, $W_4$ are all symmetric matrices, and
 \bea U_4W_4=I_{2r-\lm},\lb{d36}\\
      D_1D_2^T=I_{2p-2r},\lb{d37}\\
      U_4\td{E}_2=G_2^TD_1,\lb{d38}\eea
So
 \bea Q_4=\left(\begin{array}{cccccc}
   U_1&0&0&I_\lm&0&0\\0&U_4&\td{F}_2&0&I_{2r-\lm}&0\\0&0&D_1&0&0&0\\\Lm&0&0&W_1&0&0
   \\0&0&\td{E}_2&0&W_4&0\\0&G_2&\td{D}_3&0&K_2&D_2\end{array}\right).\eea

By (\ref{d36})-(\ref{d38}), we have both $\td{E}_2$ and $G_2$ are
zero or nonzero. By definition 2.3 we have $Q_4\sim Q_3\sim Q$. Then
by (\ref{d35}),
$\left(\begin{array}{ccc}\Lm&0&0\\0&0&\td{E}_2\\0&G_2&\td{D}_3\end{array}\right)$
is invertible. So both $\td{E}_2$ and $G_2$ are nonzero.

Since $Q_4$ is symplectic, by  (\ref{d38}) we have
 \be
 \left(\begin{array}{ccc}U_1&0&0\\0&U_4&\td{F}_2\\0&0&D_1\end{array}\right)^T
 \left(\begin{array}{ccc}\Lm&0&0\\0&0&\td{E}_2\\0&G_2&\td{D}_3\end{array}\right)
  =\left(\begin{array}{ccc}U_1\Lm&0&0\\0&0&U_4\td{E}_2\\0&(U_4\td{E}_2)^T&D_1^T\td{D}_3+\td{B}_2^T\td{E}_2
  \end{array}\right)\lb{d39}\ee
which is a symmetric matrix.

Denote by
$F=\left(\begin{array}{cc}0&U_4\td{E}_2\\(U_4\td{E_2})^T&D_1^T\td{D}_3+\td{B}_2^T\td{E}_2\end{array}\right)$.
Since $U_4\td{E}_2$ is nonzero, in the following we prove that
$m^+(F)\ge 1$.

Note that here $U_4\td{E}_2$ is a $(2r-\lm)\times (2p-2r)$ matrix
and $D_1^T\td{D}_3+\td{B}_2^T\td{E}_2$ is a $(2p-2r)\times (2p-2r)$
matrix. Denote by $U_4\td{E}_2=(e_{ij})$ and
$D_1^T\td{D}_3+\td{B}_2^T\td{E}_2=(d_{ij})$, where $e_{ij}$ and
$d_{ij}$ are elements on the $i$-th row and $j$-th column of the
corresponding matrix. Since $U_4\td{E}_2$ is nonzero, there exist an
$e_{ij}\neq 0$ for some $1\le i\le 2r-\lm$ and $1\le j\le 2p-2r$.
Let $x= (0,..,0,e_{ij},0,...0)^T\in \R^{2r-\lm}$ whose $i$-th row is
$ e_{ij}$ and other rows are all zero, and
$y=(0,...,0,\rho,0,...,0)^T\in\R^{2p-2r}$ whose $j$-th row is $\rho$
and other rows are all zero. Then we have
 \bea F\left(\begin{array}{c}x\\y\end{array}\right)\cdot \left(\begin{array}{c}x\\y\end{array}\right)
    = 2\rho e_{ij}^2-\rho^2d_{jj}>0\nn\eea
for $\rho>0$ is small enough. Hence the dimension of positive
definite space of $F$ is at least 1, thus $m^+(F)\ge 1$.
Then
 \be m^+\left(\left(\begin{array}{ccc}U_1\Lm&0&0\\0&0&U_4\td{E}_2\\0&(U_4\td{E}_2)^T&D_1^T\td{D}_3+\td{B}_2^T\td{E}_2
     \end{array}\right)\right)=m^+(\Lm)+m^+(F)\ge 1.\lb{d40}\ee
Then by (\ref{d39}), (\ref{d40}) and (ii) of Lemma 2.5, we have
 \be \frac{1}{2}\sgn M_\var(Q_4)\le 2p-1=n-2,\quad 0<-\var\ll 1.\lb{d41}\ee
Since $Q\sim Q_4$, by (\ref{d41}) and Lemma 2.4 we have
 \be \frac{1}{2}\sgn M_\var(Q)\le 2p-1, 0<-\var\ll 1.\ee
 Then since $P\sim (-I_2)\diamond Q$, by Theorem 2.1, Remark 2.2 and Lemma 2.4 we have
 \bea
  &&(i_{L_0}(\ga)+\nu_{L_0}(\ga))-(i_{L_1}(\ga)+\nu_{L_1}(\ga))\nn\\
  &=& \frac{1}{2} M_\var(P)\nn\\
  &=& \frac{1}{2}\sgn M_\var((-I_2)\diamond Q)\nn\\
  &=& \frac{1}{2}\sgn M_\var(-I_2)+\frac{1}{2}\sgn M_\var(Q)\nn\\
  &=&0+\frac{1}{2}\sgn M_\var(Q)\nn\\
  &\le& n-2.\lb{d42}\eea
Thus (\ref{d9}) and (\ref{d42}) yields a contradiction. And in Case
II we can always obtain a contradiction.

{\bf Case III.} $i(\ga^2)=n+2k+1$, where $0\le k\le p-1$.

 In this case by (\ref{d2}) we have
   \be i_{L_0}(\ga)+i_{L_1}(\ga)=2k+1.\ee
   Since $i_{L_0}(\ga)\ge 0$ and $i_{L_1}(\ga)\ge 0$  we can write $i_{L_0}(\ga)=k+1+r$ and $i_{L_1}(\ga)=k-r$ for
   some integer $-k\le r \le k$. Then by (\ref{d9}) we have
    \be n-1\ge \nu_{L_0}(\ga)-\nu_{L_1}(\ga)=n-2r-2.\lb{d45}\ee
    Thus $r\ge 0$ and $0\le r\le k$.

 By Theorem 2.1 and (i)
  of Lemma 2.5 we have
   \be 2r+1=i_{L_0}(\ga)-i_{L_1}(\ga)=\frac{1}{2}M_\var (P)\le
   n-\nu_{L_0}(\ga)\ee
   which yields $\nu_{L_0}(\ga)\le n-2r-1$. Then by (\ref{d45}) and $\nu_{L_1}(\ga)\ge 1$ we have
      \be \nu_{L_0}(\ga)=n-2r-1,\quad \nu_{L_1}(\ga)=1.\lb{d46}\ee
   Then by (\ref{d8}) we have
    \be
    S_{P^2}^+(1)=(n-2r-1)+\frac{1-n}{2}-(k-r)=\frac{1+n}{2}-k-r -1 =p-k-r\ge
    1.\lb{d47}\ee
  Then by (\ref{d10}) and $\nu(P^2)=\nu_{L_0}(\ga)+\nu_{L_1}(\ga)=n-2r$ and Lemma 2.1 we have
  \bea P^2\approx I_2^{\diamond (p-k-r)}\diamond
  N_1(1,-1)^{\diamond (2k+1)}\diamond
     R(\theta_1)\diamond \cdots\diamond R(\theta_q),\nn\eea
     where $q=n-(p-k-r)-(2k+1)=p+r-k\ge p-k\ge 1$.

Since in this case $\rank B=2r+1\le n-2$, by the same argument of
(iii) in  Case II, we have
  \bea Q\sim Q_1=\left(\begin{array}{cccc}
   A_1&B_1&I_{2r+1}&0\\0&D_1&0&0\\A_3&B_3&A_2&0\\C_3&D_3&C_2&D_2\end{array}\right).\nn\eea
 Then by the same argument of Subcases 1, 2, 3 of Case II, we can always obtain a
 contradiction in Case III.
The proof of Theorem 3.1 is complete.          \hfill\hb
\\
Now we are ready to give a proof of Theorem 1.1. For $\Sg\in
\mathcal{H}_b^{s,c}(2n)$, let $j_\Sg: \Sg \rightarrow[0,+\infty)$ be
the gauge function of $\Sg$ defined by
        \bea
        j_{\Sg}(0)=0,\quad {\rm and} \quad  j_\Sg(x)=\inf\{\lambda >0\mid
        \frac{x}{\lambda}\in C\}, \quad \forall x \in
        \R^{2n}\setminus\{0\},\nn
        \eea
where $C$ is the domain enclosed by $\Sg$.

 Define
        \bea H_\alpha(x)=(j_\Sg(x))^\alpha,\;\alpha>1,\quad
        H_\Sg(x)=H_2(x),\; \forall x \in
        \R^{2n}.\lb{7.2}
        \eea
Then $H_{\Sigma} \in C^2 (\R^{2n}\backslash \{0\},\R)\cap
    C^{1,1}(\R^{2n},\R)$.

We consider the following fixed energy problem \bea
\dot{x}(t) &=& JH_\Sg'(x(t)), \lb{7.4}\\ H_\Sg(x(t)) &=& 1,   \lb{7.5}\\
x(-t) &=& Nx(t),  \lb{7.6}\\ x(\tau+t) &=& x(t),\quad \forall\,
t\in\R. \lb{7.7} \eea

 Denote by
$\mathcal{J}_b(\Sg,2)\;(\mathcal{J}_b(\Sg,\alpha)$ for $\alpha=2$ in
(\ref{7.2})) the set of all solutions $(\tau,x)$ of problem
(\ref{7.4})-(\ref{7.7}) and by $\tilde{\mathcal{J}}_b(\Sg,2)$ the
set of all geometrically distinct solutions of
(\ref{7.4})-(\ref{7.7}). By Remark 1.2 of \cite{LiuZhang} or
discussion in \cite{LZZ}, elements in $\mathcal{J}_b(\Sg)$ and
$\mathcal{J}_b(\Sg,2)$ are one to one correspondent. So we have
$^\#\td{\mathcal{J}}_b(\Sg)$=$^\#\td{\mathcal{J}}_b(\Sg,2)$.

For readers' convenience in the following we list some known results
which will be used in the proof of Theorem 1.1. In the following of
this paper, we write
$(i_{L_0}(\gamma,k),\nu_{L_0}(\gamma,k))=(i_{L_0}(\gamma^k),\nu_{L_0}(\gamma^k))$
for any symplectic path $\gamma\in \mathcal {P}_{
 {\tau}}(2n)$ and $k\in \N$, where $\ga^k$ is defined by Definition 2.5. We have

\noindent{\bf Lemma 3.1.} (Theorem 1.5 and of \cite{LiuZhang} and
Theorem 4.3 of \cite{LZ}) {\it Let $\ga_j\in \mathcal {P}_{
 {\tau_j}}(2n)$ for $j=1,\cdots,q$.  Let
 $M_j=\ga^2_j(2\tau_j)=N\ga_j(\tau_j)^{-1}N\ga_j(\tau_j)$, for $j=1,\cdots,q$. Suppose
          \bea \hat{i}_{L_0}(\ga_j)>0, \quad
          j=1,\cdots,q.\nn\eea
  Then there exist infinitely many $(R, m_1, m_2,\cdots,m_q)\in \N^{q+1}$ such that

  (i) $\nu_{L_0}(\ga_j, 2m_j\pm 1)=\nu_{L_0}(\ga_j)$,

  (ii) $i_{L_0}(\ga_j, 2m_j-1)+\nu_{L_0}(\ga_j,2m_j-1)=R-(i_{L_1}(\ga_j)+n+S_{M_j}^+(1)-\nu_{L_0}(\ga_j))$,

  (iii) $i_{L_0}(\ga_j,2m_j+1)=R+i_{L_0}(\ga_j)$.

\noindent and
    (iv) $\nu(\ga_j^2, 2m_j\pm 1)=\nu(\ga_j^2)$,

  (v) $i(\ga_j^2, 2m_j-1)+\nu(\ga_j^2,
  2m_j-1)=2R-(i(\ga_j^2)+2S_{M_j}^+(1)-\nu(\ga_j^2))$,

  (vi) $i(\ga_j^2,2m_j+1)=2R+i(\ga_j^2)$,

\noindent where we have set $i(\ga_j^2, n_j)=i(\ga_j^{2n_j}, [0,
2n_j\tau_j])$, $\nu(\ga_j^2, n_j)=\nu(\ga_j^{2n_j}, [0,
2n_j\tau_j])$ for $n_j\in\N$.}\\

 \noindent{\bf Lemma 3.2}
(Lemma 1.1 of \cite{LiuZhang}) {\it Let $(\tau ,x)\in
\mathcal{J}_b(\Sg,2)$ be symmetric in the sense that
$x(t+\frac{\tau}{2})=-x(t)$ for all $t\in \R$ and $\ga$ be the
associated symplectic path of $(\tau,x)$. Set
$M=\ga(\frac{\tau}{2})$. Then there is a continuous symplectic path
      \bea
      \Psi(s)=P(s) M P(s)^{-1}, \quad s\in [0,1]\nn\eea
such that
      \bea \Psi(0)=M,\qquad \Psi(1)=(-I_2)\diamond \tilde{M},\;\;\;\; \td{M}\in
      \Sp(2n-2),\nn\eea
      \bea \nu_1(\Psi(s))=\nu_1(M), \quad \nu_2(\Psi(s))=\nu_2(M),
      \quad  \forall \;s\in [0,1],\nn\eea
 where
$P(s)=\left(\begin{array}{cc}\psi(s)^{-1}&0\\0&\psi(s)^T\end{array}\right)$
and $\psi$ is a continuous $n\times n$ matrix path with
$\det\psi(s)>0$ for all $s\in [0,1]$.}

For any $(\tau,x)\in\mathcal{J}_b(\Sg,2)$ and $m\in\N$, as in
\cite{LiuZhang} we denote by $i_{L_j}(x,m)=i_{L_j}(\ga_x^m,
[0,\frac{m\tau}{2}])$ and $\nu_{L_j}(x,m)=\nu_{L_j}(\ga_x^m,
[0,\frac{m\tau}{2}])$ for $j=0,1$ respectively.  Also we denote by
$i(x,m)=i(\ga_x^{2m}, [0,m\tau])$ and $\nu(x,m)=\nu(\ga_x^{2m},
[0,m\tau])$. If $m=1$, we denote by $i(x)=i(x,1)$ and
$\nu(x)=\nu(x,1)$. By Lemma 6.3 of \cite{LiuZhang} we have\\

\noindent {\bf Lemma 3.3.} {\it Suppose
$^\#\tilde{\mathcal{J}}_b(\Sg)<+\infty$. Then there exist an integer
$K\ge 0$ and an injection map $\phi: \N+K\mapsto
\mathcal{J}_{b}(\Sg,2)\times \N$ such that

(i) For any $k\in \N+K$, $[(\tau,x)]\in \mathcal{J}_{b}(\Sg,2)$ and
$m\in \N$ satisfying $\phi(k)=([(\tau \;,x)],m)$, there holds
              $$i_{L_0}(x,m)\le k-1\le i_{L_0}(x,m)+\nu_{L_0}(x,m)-1,$$
where $x$ has minimal period $\tau$.

(ii) For any $k_j\in \N+K$, $k_1<k_2$, $(\tau_j,x_j)\in \mathcal
{J}_b(\Sg,2)$ satisfying $\phi(k_j)=([(\tau_j \;,x_j)],m_j)$ with
$j=1,2$ and $[(\tau_1 \;,x_1)]=[(\tau_2 \;,x_2)]$, there holds
       $$m_1<m_2.$$}

\noindent{\bf Lemma 3.4.} (Lemma 7.2 of \cite{LiuZhang}) {\it Let
$\ga\in \P_\tau(2n)$ be extended to $[0,+\infty)$ by
$\ga(\tau+t)=\ga(t)\ga(\tau)$ for all $t>0$. Suppose
$\ga(\tau)=M=P^{-1}(I_2\diamond \td{M})P$ with $\td{M}\in \Sp(2n-2)$
and $i(\ga)\ge n$. Then we have
      \bea i(\ga,2)+2S_{M^2}^+(1)-\nu(\ga,2)\ge n+2.\nn\eea}

\noindent{\bf Lemma 3.5} (Lemma 7.3 of \cite{LiuZhang}) {\it For any
$(\tau,x)\in \mathcal{J}_b(\Sg,2)$ and $m\in \N$, we have
\bea i_{L_0}(x,m+1)-i_{L_0}(x,m)&\ge& 1,\nn\\
    i_{L_0}(x,m+1)+\nu_{L_0}(x,m+1)-1&\ge&
   i_{L_0}(x,m+1)>i_{L_0}(x,m)+\nu_{L_0}(x,m)-1.\nn\eea}

{\bf Proof of Theorem 1.1.} By Theorem 1.1 of \cite{LiuZhang} we
have $^{\#}\td{\J}_b(\Sg)\ge \left[\frac{n}{2}\right]+1$ for
$n\in\N$. So we only need to prove Theorem q.q for the case $n\ge 3$
and $n$ is odd. The method of the proof is similar as that of
\cite{LiuZhang}.

It is suffices to consider the case
 $^\#\tilde{\mathcal{J}}_b(\Sg)<+\infty$. Since $-\Sg=\Sg$, for
 $(\tau,x) \in \mathcal{J}_b(\Sg,2)$  we have
          \bea &&H_\Sg(x)=H_\Sg(-x),\nn\\
               &&H_\Sg'(x)=- H_\Sg'(-x),\nn\\
                &&H_\Sg''(x)= H_\Sg''(-x).\lb{8.43}\eea
So $(\tau,-x)\in \mathcal{J}_b(\Sg,2)$. By (\ref{8.43}) and the
definition of $\ga_x$ we have that
 \bea \ga_x=\ga_{-x}.\nn\eea
So we have
       \bea &&(i_{L_0}(x,m),\nu_{L_0}(x,m))=(i_{L_0}(-x,m),\nu_{L_0}(-x,m)),\nn\\
       &&(i_{L_1}(x,m),\nu_{L_1}(x,m))=(i_{L_1}(-x,m),\nu_{L_1}(-x,m)),\quad \forall m\in
       \N.\lb{8.45}\eea
 So we can write
  \be \td{\mathcal {J}}_b(\Sg,2)=\{[(\tau_j,x_j)]|
j=1,\cdots,p\}\cup\{[(\tau_k,x_k)],[(\tau_k,-x_k)]|k=p+1,\cdots,p+q\}.\lb{8.46}\ee
with $x_j(\R)=-x_j(\R)$ for $j=1,\cdots,p$ and $x_k(\R)\neq
-x_k(\R)$ for $k=p+1,\cdots,p+q$. Here we remind that $(\tau_j,x_j)$
has minimal period $\tau_j$ for $j=1,\cdots,p+q$ and
$x_j(\frac{\tau_j}{2}+t)=-x_j(t), \;t\in\R$ for $j=1,\cdots,p$.

 By Lemma 3.3 we have an integer $K\ge 0$ and an injection map
 $\phi: \N+K\to
\mathcal{J}_b(\Sg,2)\times \N$. By (\ref{8.45}), $(\tau_k,x_k)$ and
$(\tau_k,-x_k)$ have the same $(i_{L_0},\nu_{L_0})$-indices.
 So by Lemma 3.3,
 without loss of generality, we can further require that
       \bea {\rm Im} (\phi)\subseteq \{[(\tau_k,x_k)]|k=1,2,\cdots,p+q\}\times
       \N.\lb{8.47}\eea
By the strict convexity of $H_\Sg$ and (6.19) of \cite{LiuZhang}),
we have
        \bea \hat{i}_{L_0}(x_k)>0,\quad
        k=1,2,\cdots,p+q.\nn\eea
Applying  Lemma 3.1 to the following associated symplectic paths
$$\ga_1,\;\cdots,\;\ga_{p+q},\; \ga_{p+q+1},\;\cdots,\;\ga_{p+2q}$$
of
$(\tau_1,x_1),\;\cdots,\;(\tau_{p+q},x_{p+q}),\;(2\tau_{p+1},x_{p+1}^2),\;\cdots,\;
     (2\tau_{p+q},x_{p+q}^2)$ respectively,
there exists a vector $(R,m_1,\cdots,m_{p+2q})\in \N^{p+2q+1}$ such
that $R>K+n$ and
   \bea &&i_{L_0}(x_k, 2m_k+1)=R+i_{L_0}(x_k),\lb{8.49}\\
        && i_{L_0}(x_k,2m_k-1)+\nu_{L_0}(x_k,2m_k-1)\nn\\
        &=&R-(i_{L_1}(x_k)+n+S_{M_k}^+(1)-\nu_{L_0}(x_k)),\lb{8.50}\eea
    for $k=1,\cdots,p+q,$ $M_k=\ga_k^2(\tau_k)$, and
     \bea &&i_{L_0}(x_k, 4m_k+2)=R+i_{L_0}(x_k,2),\lb{8.51}\\
         &&i_{L_0}(x_k,4m_k-2)+\nu_{L_0}(x_k,4m_k-2)\nn\\
         &=&R-(i_{L_1}(x_k,2)+n+S_{M_k}^+(1)-\nu_{L_0}(x_k,2)),\lb{8.52}\eea
      for $k=p+q+1,\cdots,p+2q$ and $M_k=\ga_k^4(2\tau_k)=\ga_k^2(\tau_k)^2$.

By Lemma 3.1, we also have \bea i(x_k,
2m_k+1)&=&2R+i(x_k),\lb{8.53}\\
         i(x_k,2m_k-1)+\nu(x_k,2m_k-1)
        &=&2R-(i(x_k)+2S_{M_k}^+(1)-\nu(x_k)),\lb{8.54}\eea
    for $k=1,\cdots,p+q,$ $M_k=\ga_k^2(\tau_k)$, and
     \bea i(x_k, 4m_k+2)&=&2R+i(x_k,2),\lb{8.55}\\
        i(x_k,4m_k-2)+\nu(x_k,4m_k-2)
         &=&2R-(i(x_k,2)+2S_{M_k}^+(1)-\nu(x_k,2)),\lb{8.56}\eea
      for $k=p+q+1,\cdots,p+2q$ and $M_k=\ga_k^4(2\tau_k)=\ga_k^2(\tau_k)^2$.

From (\ref{8.47}), we can set
 \bea \phi(R-(s-1))=([(\tau_{k(s)}, x_{k(s)})],m(s)),\qquad
    \forall s\in S:=\left\{1,2,\cdots,\left[\frac{n+1}{2}\right]+1\right\},\nn\eea
where $k(s)\in \{1,2,\cdots,p+q\}$ and $m(s)\in \N$.

We continue our proof  to study the symmetric and asymmetric orbits
separately. Let \bea S_1=\{s\in S|k(s)\le p\},\qquad S_2=S\setminus
S_1.\nn\eea
 We shall prove that
$^\#S_1\le p$ and $^\#S_2\le 2q$, together with the definitions of
$S_1$ and $S_2$, these yield Theorem 1.1.

\noindent{\bf Claim 1.} $^\#S_1\le p$.

\noindent {\it Proof of Claim 1.} By the definition of $S_1$,
$([(\tau_{k(s)},
 x_{k(s)})],m(s))$ is symmetric when $k(s)\le p$. We further prove
 that $m(s)=2m_{k(s)}$ for $s\in S_1$.

  In fact, by the definition of $\phi$ and Lemma 3.3,
   for all $s=1,2,\cdots,\left[\frac{n+1}{2}\right]+1$ we have
         \bea  i_{L_0}(x_{k(s)},m(s))&\le & (R-(s-1))-1=R-s \nn\\
         &\le &
         i_{L_0}(x_{k(s)},m(s))+\nu_{L_0}(x_{k(s)},m(s))-1.\lb{8.59}\eea
 By the strict convexity of $H_\Sg$ and Lemma 2.2, we have $i_{L_0}(x_{k(s)})\ge 0$, so there holds
   \bea i_{L_0}(x_{k(s)},m(s))\le R-s< R\le R+i_{L_0}(x_{k(s)})=i_{L_0}(x_{k(s)},2m_{k(s)}+1),\lb{8.60}\eea
 for every $s=1,2,\cdots,\left[\frac{n+1}{2}\right]+1$, where we have used
 (\ref{8.49}) in the last equality. Note that the proofs of (\ref{8.59}) and
 (\ref{8.60}) do not depend on the condition $s\in S_1$.

By Lemma 3.2, $\ga_{x_k}$ satisfies conditions of Theorem 3.1 with
$\tau=\frac{\tau_k}{2}$. Note that by definition
$i_{L_1}(x_k)=i_{L_1}(\ga_{x_k})$ and
$\nu_{L_0}(x_k)=\nu_{L_0}(\ga_{x_k})$. So by Theorem 3.1 we have
   \be
i_{L_1}(x_k)+S_{M_k}^+(1)-\nu_{L_0}(x_k)>\frac{1-n}{2},\quad \forall
k=1,\cdots,p.\lb{8.64}\ee
 Also for $1\le
s\le \left[\frac{n+1}{2}\right]+1$, we have
   \be -\frac{n+3}{2}= -\left(\left[\frac{n+1}{2}\right]+1\right)\le
   -s.\lb{8.65}\ee
Hence by (\ref{8.59}),(\ref{8.64}) and(\ref{8.65}), if $k(s)\le p$
we have
 \bea
&&i_{L_0}(x_{k(s)},2m_{k(s)}-1)+\nu_{L_0}(x_{k(s)},2m_{k(s)}-1)-1\nn\\
&=&
R-(i_{L_1}(x_{k(s)})+n+S_{M_{k(s)}}^+(1)-\nu_{L_0}(x_{k(s)}))-1\nn\\
&<&R-\frac{1-n}{2}-1-n=R-\frac{n+3}{2}\le R-s\nn\\
&\le&
i_{L_0}(x_{k(s)},m(s))+\nu_{L_0}(x_{k(s)},m(s))-1.\lb{8.66}\eea
 Thus
by (\ref{8.60}) and (\ref{8.66}) and Lemma 3.5 of \cite{LiuZhang} we
have
 \be 2m_{k(s)}-1< m(s)<2m_{k(s)}+1.\lb{8.67}\ee
 Hence
 \be m(s)=2m_{k(s)}.\lb{8.68}\ee
So we have
 \be \phi(R-s+1)=([(\tau_{k(s)},x_{k(s)})],2m_{k(s)}),\qquad \forall
 s\in S_1.\lb{8.69}\ee
Then by the injectivity of $\phi$, it induces another injection map
 \be \phi_1:S_1\rightarrow \{1,\cdots,p\}, \;s\mapsto k(s).\lb{8.70}\ee
 There for $^\#S_1\le p$. Claim 1 is proved.

 \noindent{\bf Claim 2.} $^\#S_2\le 2q$.

\noindent{\it Proof of Claim 2.} By the  formulas
(\ref{8.53})-(\ref{8.56}), and (59) of \cite{LLZ} (also Claim 4 on
p. 352 of \cite{Long1}), we have \be m_k=2m_{k+q}\quad {\rm for}\;\;
k=p+1,p+2,\cdots,p+q.\lb{8.71}\ee
 We set $\mathcal {A}_k=i_{L_1}(x_k,2)+S_{M_k}^+(1)-\nu_{L_0}(x_k,2)$
and $\mathcal {B}_k=i_{L_0}(x_k,2)+S_{M_k}^+(1)-\nu_{L_1}(x_k,2)$,
$p+1\le k\le p+q$, where $M_k=\ga_k(2\tau_k)=\ga(\tau_k)^2$. By
(\ref{d2}), we have
 \be
\mathcal {A}_k+\mathcal
{B}_k=i(x_k,2)+2S_{M_k}^+(1)-\nu(x_k,2)-n,\;\;\;p+1\le k\le p+q
.\lb{8.72}\ee By similar discussion of the proof of Lemma 3.2, for
any $p+1\le k\le p+q$ there exist $P_k\in \Sp(2n)$ and $\td{M}_k\in
\Sp(2n-2)$ such that \bea \ga(\tau_k)=P_k^{-1}(I_2\diamond
\td{M}_k)P_k.\nn\eea
 Hence by Lemma 3.4 and (\ref{8.72}), we
have \be \mathcal {A}_k+\mathcal {B}_k\ge n+2-n=2.\lb{8.74}\ee By
Theorem 2.1, there holds
 \bea
|\mathcal {A}_k-\mathcal
{B}_k|&=&|(i_{L_0}(x_k,2)+\nu_{L_0}(x_k,2))-(i_{L_1}(x_k,2)+\nu_{L_1}(x_k,2))|\le
n.\lb{8.75}\eea
 So by (\ref{8.74}) and (\ref{8.75}) we have \be
\mathcal {A}_k\ge \frac{1}{2}((\mathcal {A}_k+\mathcal
{B}_k)-|\mathcal {A}_k-\mathcal {B}_k|)\ge \frac{2-n}{2},\quad
p+1\le k\le p+q.\lb{8.76}\ee
 By
(\ref{8.52}), (\ref{8.59}), (\ref{8.65}), (\ref{8.71}) and
(\ref{8.76}), for $p+1\le k(s)\le p+q$ we have
\bea &&i_{L_0}(x_{k(s)},2m_{k(s)}-2)+\nu_{L_0}(x_{k(s)},2m_{k(s)}-2)-1\nn\\
     &=&i_{L_0}(x_{k(s)},4m_{k(s)+q}-2)+\nu_{L_0}(x_{k(s)},4m_{k(s)+q}-2)-1\nn\\
     &=&R-(i_{L_1}(x_{k(s)},2)+n+S_{M_{k(s)}}^+(1)-\nu_{L_0}(x_{k(s)},2))-1\nn\\
     &=&R-\mathcal{A}_{k(s)}-1-n\nn\\
     &\le&R- \frac{2-n}{2}-1-n\nn\\
     &=& R-(2+\frac{n}{2})\nn\\
     &<& R-\frac{n+3}{2}\nn\\&\le& R-s\nn\\
     &\le&
     i_{L_0}(x_{k(s)},m(s))+\nu_{L_0}(x_{k(s)},m(s))-1.\lb{8.77}\eea
Thus by (\ref{8.60}), (\ref{8.77}) and Lemma 3.5, we have
     \bea 2m_{k(s)}-2<m(s)<2m_{k(s)}+1,\qquad p<k(s)\le
     p+q.\nn\eea
So \bea m(s)\in \{2m_{k(s)}-1,2m_{k(s)}\}, \qquad {\rm
for}\;\;p<k(s)\le p+q.\}\nn\eea
 Especially this yields that for any $s_0$ and $s\in
S_2$, if $k(s)=k(s_0)$, then
 \bea m(s)\in
\{2m_{k(s)}-1,2m_{k(s)}\}=\{2m_{k(s_0)}-1,2m_{k(s_0)}\}.\nn\eea Thus
by the injectivity of the map $\phi$ from Lemma 3.3, we have \bea
^\#\{s\in S_2|k(s)=k(s_0)\}\le 2\nn\eea which yields Claim 2.

     By Claim 1 and Claim 2, we have
     \bea
     ^\#\td{\mathcal{J}}_b(\Sg)=^\#\td{\mathcal{J}}_b(\Sg,2)=p+2q\ge
     ^\#S_1+^\#S_2
     =\left[\frac{n+1}{2}\right]+1.\nn\eea
The proof of Theorem 1.1 is complete. \hfill\hb

{\bf Proof of Theorem 1.2.} By \cite{LLZ}, there are at least $n$
closed characteristics on every $C^2$ compact convex central
symmetric hypersurface $\Sg$ of $\R^{2n}$. Hence by Example 1.1 the
assumption of Theorem 1.2 is reasonable. Here we prove the case
$n=5$, the proof of the case $n=4$ is the same.

We call a closed characteristic $x$ on $\Sg$ a {\it dual brake
orbit} on $\Sg$ if $x(-t)=-Nx(t)$. Then by the similar proof of
Lemma 3.1 of \cite{Zhang2}, a closed characteristic $x$ on $\Sg$ can
became a dual brake orbit after suitable time translation if and
only if $x(\R)=-Nx(\R)$. So by Lemma 3.1 of \cite{Zhang2} again, if
a closed characteristic $x$ on $\Sg$ can both became brake orbits
and dual brake orbits after suitable translation, then
$x(\R)=Nx(\R)=-Nx(\R)$, Thus $x(\R)=-x(\R)$.

Since we also have $-N\Sg=\Sg$, $(-N)^2=I_{2n}$ and $(-N)J=-J(-N)$,
dually by the same proof of Theorem 1.1, there are at least
$[(n+1)/2]+1=4$ geometrically distinct dual brake orbits on $\Sg$.

If there are exactly 5 closed characteristics on $\Sg$. By Theorem
1.1, four  closed characteristics of them must be brake orbits after
suitable time translation, then the fifth, say $y$, must be brake
orbits after suitable time translation, otherwise $Ny(-\cdot)$ will
be the sixth geometrically distinct closed characteristic on $\Sg$
which yields a contradiction. Hence all closed characteristics on
$\Sg$ must be brake orbits on $\Sg$. By the same argument we can
prove that all closed characteristics on $\Sg$ must be dual brake
orbits on $\Sg$. Then by the argument in the second paragraph of the
proof of this theorem, all these five closed characteristics on
$\Sg$ must be symmetric. Hence all of them bust be symmetric brake
orbits after suitable time translation. Thus we have proved the case
$n=5$ of Theorem 1.2 and the proof of Theorem 1.2 is complete.
\hfill\hb

\bibliographystyle{abbrv}

\end{document}